\documentclass[a4paper,11pt,reqno]{amsart}
\usepackage{amsmath,amssymb,amsthm,enumerate}
\usepackage{algorithm,algpseudocode}
\usepackage[bookmarks,colorlinks, citecolor = blue]{hyperref}
\usepackage[x11names,rgb]{xcolor}
\usepackage{caption}
\usepackage{subfig}
\usepackage{tikz}
\usetikzlibrary{arrows,shapes,matrix,decorations.pathmorphing}

\newcommand{\PP}{\mathbb{P}}
\newcommand{\K}{\mathbb{K}}
\newcommand{\hilb}{{\mathcal{H}\textnormal{ilb}}}
\newcommand{\h}{{\mathit{hilb}}}
\newcommand{\hilbp}{{\mathcal{H}\textnormal{ilb}_{p(t)}^n}}
\newcommand{\hilbd}{{\mathcal{H}\textnormal{ilb}_{d}^n}}
\newcommand{\Proj}{\textnormal{Proj}\,}
\newcommand{\up}{\textrm{e}^{+}}
\newcommand{\down}{\textrm{e}^{-}}
\newcommand{\poset}{\mathcal{P}(n,r)}

\newcommand{\NI}{\mathcal{N}(I)}

\newtheorem{theorem}{Theorem}[section]
\newtheorem{corollary}[theorem]{Corollary}
\newtheorem{proposition}[theorem]{Proposition}
\newtheorem{lemma}[theorem]{Lemma}

\theoremstyle{definition}
\newtheorem{definition}[theorem]{Definition}
\newtheorem{notation}[theorem]{Notation}
\newtheorem{example}[theorem]{Example}
\newtheorem{remark}[theorem]{Remark}

\begin{document}
 
\title{A network of rational curves on the Hilbert scheme}

\author[P.~Lella]{Paolo Lella}
\address{Universit\`a degli Studi di Torino \\ Dipartimento di Matematica \\ Via Carlo Alberto 10 \\ 10123 Torino \\ Italy}
\email{\href{mailto:paolo.lella@unito.it}{paolo.lella@unito.it}}
\urladdr{\url{http://www.dm.unito.it/dottorato/dottorandi/lella/}}


\begin{abstract}
In this paper we introduce an effective method to construct rational deformations between couples of Borel-fixed ideals. These deformations are governed by flat families, so that they correspond to rational curves on the Hilbert scheme. Looking globally at all the deformations among Borel-fixed ideals defining points on the same Hilbert scheme, we are able to give a new proof of the connectedness of the Hilbert scheme and to introduce a new criterion to establish whenever a set of points defined by Borel ideals lies on a common component of the Hilbert scheme. The paper contains a detailed algorithmic description of the technique and all the algorithms are made available.
\end{abstract}

\maketitle

\section*{Introduction}
Borel-fixed ideals are special ideals hard studied for their strong and advantageous combinatorial properties. Moreover in the study of the Hilbert scheme $\hilbp$, Borel ideals play a key role because of the result of Galligo (see \cite{Green}) saying that any component and any intersection of components of the Hilbert scheme contains a point corresponding to a subscheme defined by a Borel ideal. For instance, they are the base object in the proof by Hartshorne \cite{H66} of the connectedness of the Hilbert scheme and more recently they have been used as starting point to compute an open covering of $\hilbp$ \cite{LR,R}.

The aim of this paper is to construct a network of deformations between couples of Borel ideals. Hartshorne \cite{H66} and Peeva-Stillman \cite{PS} introduced deformations in order to move from any Borel-fixed ideal toward the lexicographic ideal: Hartshorne through \emph{fans} and Peeva and Stillman through Gr\"obner deformations. Our point of view is slightly different because initially the \lq\lq direction\rq\rq\ of the deformation does not matter. We look for a criterion for swapping couples of monomials, one inside the ideal and one outside, being careful to preserve the Borel condition and the Hilbert polynomial. In Section \ref{sec:0case} we expose a method to detect \lq\lq good\rq\rq couples of monomial for the zero-dimensional ideal and then we extend it to the general case in Section \ref{sec:genCase}. The technique introduced is completely combinatorial and avoids other tools as Gr\"obner computations.

After having found this criterion, we prove that there exists a rational deformation governing the exchange of monomials (Theorem \ref{th:flatZeroDim} and Theorem \ref{th:flatGeneralDim}) that corresponds to a rational curve on the Hilbert scheme Theorem \ref{th:curveGeneralDim}).

In Section \ref{sec:connected}, we slightly modify our technique in order to compare it with that of Peeva and Stillman. We show an alternative proof of the connectedness of the Hilbert scheme, based on the idea of moving toward a hilb-segment ideal, a generalization (introduced in \cite{CLMR}) of the notion of lexsegment ideal. 
 
Section \ref{sec:components} is devoted to the study of the components of the Hilbert scheme and their intersections. Considering the rational deformations as a whole, we introduce a new criterion to establish whenever a set of Borel points lies on a common component.

Throughout the paper, the perspective is to project effective methods for concrete computations, so every construction introduced comes with an algorithm (and its pseudo-code description) and with many examples.

\section{General Settings}
Let $\K$ be an algebraically closed field of characteristic 0, $\K[x] = \K[x_0,\ldots,x_n]$ the polynomial ring with $n+1$ variables and let $\PP^n = \Proj \K[x]$ be the $n$-dimensional projective space.

$\preceq$ will be a graded term ordering on the monomials of $\K[x]$ with the unusual, but convenient in our case, ordering on the variables $x_n \succ \ldots \succ x_0$. Given a set of polynomials $G = \{f_1,\ldots,f_s\}$, we will denote by $\langle G\rangle$ the ideal generated by those polynomials.

To denote monomials in $\K[x]$, we will use the compact notation $x^\alpha = x_n^{\alpha_n} \cdots x_0^{\alpha_0}$, for every multi-index $\alpha = (\alpha_n,\ldots,\alpha_0)$ with non-negative entries.

\medskip
 
$\hilbp$ will be used to denote the Hilbert scheme parametrizing the family $\h$
\begin{center}
\begin{tikzpicture}
\node (a) at (0,1.5) {$\mathcal{X}$};
\node (b) at (3,1.5) {$\PP^n \times \hilbp$};
\node (c) at (3,0) {$\hilbp$};
\draw [right hook->] (a) -- (b);
\draw [->] (a) --node[below]{$\h$} (c);
\draw [->] (b) -- (c);
\end{tikzpicture}
\end{center}
of subschemes of the projective space $\PP^n$ with Hilbert polynomial equal to the polynomial $p(t)$. We will often denote by $r$ the \emph{Gotzmann number} of the Hilbert polynomial $p(t)$, that is the best upper bound for the Castelnuovo-Mumford regularity of a subscheme having such Hilbert polynomial. Let us remark that the Gotzmann number coincides with the regularity of the unique saturated lexsegment ideal having Hilbert polynomial $q(t) = \binom{n+t}{n}-p(t)$ (see \cite{Green}).

$\Delta p(t) = p(t)-p(t-1)$ will be the first difference of the polynomial $p(t)$. Let us remind that given an ideal $I$ such that $S/I$ has Hilbert polynomial $p(t)$, for a generic linear form $\ell$ we have the short exact sequence:
\begin{equation*}
 0 \ \longrightarrow\ {S}/{I}(t-1) \ \stackrel{\cdot\, \ell}{\longrightarrow}\ {S}/{I}(t) \ \longrightarrow\ {S}/{(I,\ell)}(t) \ \longrightarrow\ 0;
\end{equation*}
so $\Delta p(t)$ is the Hilbert polynomial of the generic hyperplane section of $S/I$. We assume $\Delta^0 p(t) = p(t)$ and recursively we define $\Delta^m p(t) = \Delta \big( \Delta^{m-1}p(t)\big)$. Of course $\Delta^{\deg p(t) + 1} p(t) = 0$. 

We refer to \cite{B,HS,Sern} for the classical construction of the Hilbert scheme as a subscheme of the Grassmannian $\mathbb{G}\left(q(r),\binom{n+r}{n}\right)$. However we recall the fundamental result by Gotzmann used for the explicit construction. 
\begin{theorem}[Gotzmann's Persistence Theorem]\label{th:Gotzmann}
 Let $p(t)$ be an admissible Hilbert polynomial and let $I \subseteq \K[x]$ a homogeneous ideal generated in degree $\leqslant r$ such that $\dim_k I_r = q(r)$. Then:
\begin{equation}
 \dim_k I_{r+1} = q(r+1) \qquad \Longrightarrow \qquad  \dim_k I_{t} = q(t),\ \forall\ t \geqslant r. 
\end{equation}
\end{theorem}
So, given an ideal $I = I_{\geqslant r}$, generated in degree $r$, corresponding to a point of the Grassmannian $\mathbb{G}\left(q(r),\binom{n+r}{n}\right)$, to know if $\Proj \K[x]/I$ belongs to $\hilbp$, it suffices to check the condition $\dim_k I_{r+1} \leqslant q(r+1)$,  keeping in mind the Macaulay's Estimates of Growth of Ideals \cite[Theorem 3.3]{Green} that ensures the other inequality $\dim_k I_{r+1} \geqslant q(r+1)$.

Since we will always work with ideals corresponding to subschemes parametrized by $\hilbp$, by abuse of notation we will say that an ideal $I$ belongs to $\hilbp$ meaning that the subscheme $\Proj \K[x]/I$ is in $\hilbp$ and that $I = I_{\geqslant r}$ is generated in degree $r$ equal to the Gotzmann number of $p(t)$.

\medskip

In this context, a special role is up to Borel-fixed ideals. A \emph{Borel ideal} is defined as an ideal fixed by the action of the Borel subgroup of the linear group $GL(n+1)$ (upper triangular matrices). By Galligo's Theorem \cite[Theorem 1.27]{Green}, it is well-know that every component and every intersection of components of $\hilbp$ contains at least one among these ideals. In a certain sense, we can consider them as distributed all over the Hilbert scheme.

\section{Combinatorial interpretation of Borel ideals}
In the case of a polynomial ring with coefficients in an algebraically closed field of characteristic 0 and so in our context, Borel-fixed ideals are also called \emph{strongly stable} ideals and they can be characterized by the following combinatorial property (see \cite{Green}).
\begin{definition}\label{def:BorelIdeal}
 A homogeneous ideal $I \subset \K[x]$ is Borel-fixed if
\begin{enumerate}
 \item $I$ is a monomial ideal;
 \item $\forall\ x^\alpha \in I$, $x_j \dfrac{x^\alpha}{x_i}\in I$, $\forall\ i$ s.t. $x_i \mid x^\alpha$ and $\forall\ j > i$.
\end{enumerate}
\end{definition}

Starting from this definition, we recall the \emph{Borel elementary moves} already used in \cite{Green,CLMR}.
\begin{definition}\label{def:ElementaryMoves}
Given any monomial $x^\alpha \in \K[x]$, we define
\begin{itemize}
 \item the (admissible) elementary decreasing move $\down_i$ as
\[
 \down_i (x^\alpha) = \down_i(x_n^{\alpha_n} \cdots x_i^{\alpha_i} \cdots x_0^{\alpha_0}) = x_n^{\alpha_n} \cdots x_i^{\alpha_i-1} x_{i-1}^{\alpha_{i-1}+1} \cdots x_0^{\alpha_0},\ \forall\ 0 < i \leqslant n \text{ s.t. } \alpha_i  > 0;
\]
\item the (admissible) elementary increasing move $\up_j$ as
\[
 \up_j (x^\alpha) = \up_j(x_n^{\alpha_n} \cdots x_j^{\alpha_j} \cdots x_0^{\alpha_0}) = x_n^{\alpha_n} \cdots x_{j+1}^{\alpha_{j+1}+1} x_j^{\alpha_j-1} \cdots x_0^{\alpha_0},\ \forall\ 0 \leqslant j < n \text{ s.t. } \alpha_j  > 0.
\]
\end{itemize} 

Since $\down_i (x^\alpha) = \frac{x_{i-1}}{x_i} x^\alpha$ and $\up_j (x^\beta) = \frac{x_{j+1}}{x_j} x^\beta$, the Borel moves can be viewed as elements in $\K(x)$ and a composition of elementary moves $F$ corresponds to a monomial $x^\gamma \in \K(x)$, $\gamma \in \mathbb{Z}^{n+1}$. So a composition of moves $F$ will be admissible on the monomial $x^\alpha$ if the monomial $x^\gamma\, x^\alpha$ belongs to $\K[x]$. Going by the commutativity of the product, given a composition $F = (\down_{i_s})^{\lambda_s} \cdot \ldots \cdot (\down_{i_1})^{\lambda_1}$ admissible on $x^\alpha$, we will suppose $i_1 > \ldots > i_s$ and in the same way for $G = (\up_{j_u})^{\mu_u} \cdot \ldots \cdot (\up_{j_1})^{\mu_1}$ admissible on $x^\beta$, we will suppose $j_1 < \ldots < j_u$.

We prefer to think these moves more as maps on the monomials than as simple multiplications, because of the computational perspective we have and because this is useful for a graphical representation.
\end{definition}

We denote by $\leq_B$ the partial order induced on the monomials by the transitive closure of the relations
\begin{equation}\label{eq:BorelOrder}
\up_j(x^\alpha) \geq_B x^\alpha \geq_B \down_i(x^\alpha).
\end{equation}

\begin{proposition}
Let $x^\alpha$ and $x^\beta$ be two monomials of degree $r$.
\[ 
x^\alpha \geq_B x^\beta\quad \Longleftrightarrow\quad \sigma(j) = \sum_{i=j}^n (\alpha_i - \beta_i) \geqslant 0,\ \forall\ j.
\]
\begin{proof}
$(\Rightarrow)$ $x^\alpha \geq_B x^\beta$ means $x^\beta = (\down_{i_s})^{\lambda_s} \cdot \ldots \cdot (\down_{i_1})^{\lambda_1} (x^\alpha)$, $i_1 > \ldots > i_s$. Obviously $\sigma(0) = \vert\alpha\vert-\vert\beta\vert = 0$, $\sigma(j) = 0,\ \forall\ j > i_1$ and $\sigma(i_1) > 0$ because $\alpha_{i_1} > \beta_{i_1}$. Let $x^\gamma = (\down_{i_1})^{\lambda_1}(x^\alpha)$. By construction $\gamma = (\alpha_0,\ldots,\alpha_{i_1-1} + (\alpha_{i_1}-\beta_{i_1}),\beta_{i_1},\ldots,\alpha_n)$, that is $\lambda_1 = \alpha_{i_1}-\beta_{i_1}$. Repeating the reasoning on $x^\gamma \geq_B x^\beta = (\down_{i_s})^{\lambda_s} \cdot \ldots \cdot (\down_{i_2})^{\lambda_2}(x^\gamma)$, we prove $\sigma(j) \geqslant 0, \forall\ j$. 

\noindent $(\Leftarrow)$ It is sufficient to consider the composition of decreasing moves $F = \left(\down_1\right)^{\sigma(1)}\cdot \ldots \cdot \left(\down_n\right)^{\sigma(n)}$.
\end{proof}
\end{proposition}

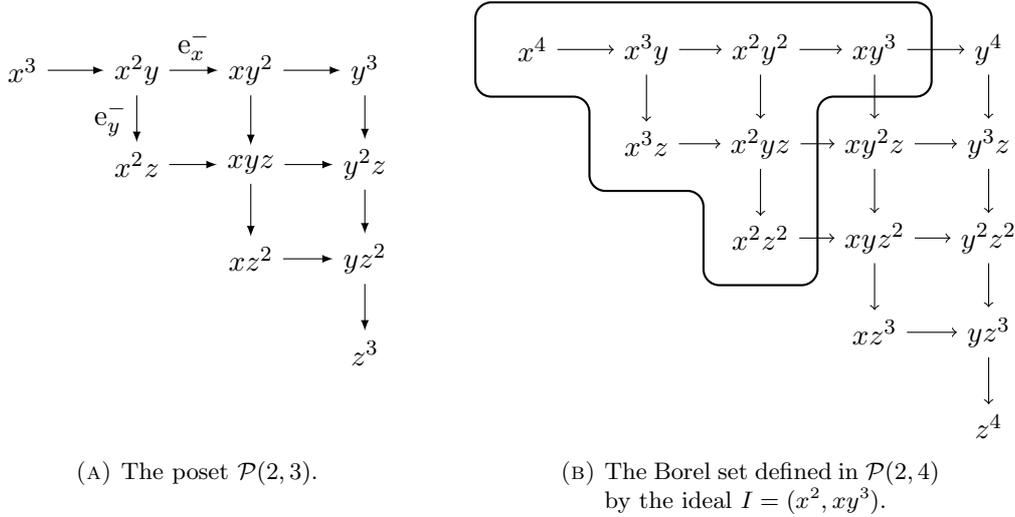
\begin{figure}[!ht] 
\begin{center}\captionsetup[subfloat]{format=hang}
\subfloat[][The poset $\mathcal{P}(2,3)$.]{
 \begin{tikzpicture}[>=latex]
\node (N_1) at (1,5.75) [] {$x^3$};
\node (N_2) at (2.5,5.75) [] {$x^2y$};
\node (N_3) at (4,5.75) [] {$xy^2 $};
\node (N_4) at (5.5,5.75) [] {$y^3$};

\node (N_5) at (2.5,4.5) [] {$x^2z$};
\node (N_6) at (4,4.5) [] {$xyz$};
\node (N_7) at (5.5,4.5) [] {$y^2z$};

\node (N_8) at (4,3.25) [] {$xz^2$};
\node (N_9) at (5.5,3.25) [] {$yz^2$};

\node (N_10) at (5.5,2) [] {$z^3$};

\node at (1,1) [] {\phantom{d}};

\draw [->] (N_1) -- (N_2);
\draw [->] (N_2) -- node[above] {$\down_x$} (N_3);
\draw [->] (N_3) -- (N_4);

\draw [->] (N_2) -- node[left] {$\down_y$}  (N_5);
\draw [->] (N_3) -- (N_6);
\draw [->] (N_4) -- (N_7);

\draw [->] (N_5) -- (N_6);
\draw [->] (N_6) -- (N_7);

\draw [->] (N_6) -- (N_8);
\draw [->] (N_7) -- (N_9);

\draw [->] (N_8) -- (N_9);

\draw [->] (N_9) -- (N_10);
\end{tikzpicture}
}
\quad\quad
\subfloat[][The Borel set defined in $\mathcal{P}(2,4)$ by the ideal $I=(x^2,xy^3)$.]{\label{fig:ExampleBorel_b}
\begin{tikzpicture}
\node (M_1) at (8,6) [] {$x^4$};
\node (M_2) at (9.5,6) [] {$x^3y$};
\node (M_3) at (11,6) [] {$x^2y^2$};
\node (M_4) at (12.5,6) [] {$xy^3$};
\node (M_5) at (14,6) [] {$y^4$};

\node (M_6) at (9.5,4.75) [] {$x^3z$};
\node (M_7) at (11,4.75) [] {$x^2yz$};
\node (M_8) at (12.5,4.75) [] {$xy^2z$};
\node (M_9) at (14,4.75) [] {$y^3z$};

\node (M_10) at (11,3.5) [] {$x^2z^2$};
\node (M_11) at (12.5,3.5) [] {$xyz^2$};
\node (M_12) at (14,3.5) [] {$y^2z^2$};

\node (M_13) at (12.5,2.25) [] {$xz^3$};
\node (M_14) at (14,2.25) [] {$yz^3$};

\node (M_15) at (14,1) [] {$z^4$};

\draw [->] (M_1) -- (M_2);
\draw [->] (M_2) -- (M_3);
\draw [->] (M_3) -- (M_4);
\draw [->] (M_4) -- (M_5);

\draw [->] (M_2) -- (M_6);
\draw [->] (M_3) -- (M_7);
\draw [->] (M_4) -- (M_8);
\draw [->] (M_5) -- (M_9);

\draw [->] (M_6) -- (M_7);
\draw [->] (M_7) -- (M_8);
\draw [->] (M_8) -- (M_9);

\draw [->] (M_7) -- (M_10);
\draw [->] (M_8) -- (M_11);
\draw [->] (M_9) -- (M_12);

\draw [->] (M_10) -- (M_11);
\draw [->] (M_11) -- (M_12);

\draw [->] (M_11) -- (M_13);
\draw [->] (M_12) -- (M_14);

\draw [->] (M_13) -- (M_14);
\draw [->] (M_14) -- (M_15);

\draw [rounded corners=6pt,thick] (7.25,6.625) -- (13.25,6.625) -- (13.25,5.375) -- (11.75,5.375) -- (11.75,2.875) -- (10.25,2.875) -- (10.25,4.125) -- (8.75,4.125) -- (8.75,5.375) -- (7.25,5.375) -- cycle;
\end{tikzpicture}
}
\end{center}
\caption{\label{fig:ExampleBorel} Examples of posets and Borel sets defined by Borel-fixed ideals.}
\end{figure}

Following the approach introduced for instance in \cite{MR99} and \cite{sherman-2005}, we will call \emph{poset}, and we will denote it by $\mathcal{P}(n,r)$, the set of monomials of degree $r$ in the polynomial ring $\K[x]$ with the order relations given by the elementary moves. Easily $\mathcal{P}(n,r)$ can be represented by a graph where the monomials correspond to the vertices and the edges represent the elementary moves (see Figure \ref{fig:ExampleBorel}).

\begin{definition}
A subset $B$ of $\mathcal{P}(n,r)$ is called \emph{Borel set} if
\begin{equation}
 \forall\ x^\alpha \in B,\quad \up_j(x^\alpha) \in B,\ \text{for every admissible elementary move } \up_j.
\end{equation}
Given a Borel set $B$, we denote by $N$ the set of monomials in $\poset \setminus B$. It holds
\begin{equation}
 \forall\ x^\beta \in N,\quad \down_i(x^\beta) \in N,\ \text{for every admissible elementary move } \down_i.
\end{equation}
\end{definition}

\begin{remark}
 Given a Borel-fixed ideal $I \subset \K[x]$, the standard monomial basis of $I_{r}$ is a Borel set for the poset $\mathcal{P}(n,r)$. In fact by Definition \ref{def:BorelIdeal}, we know that if $x^\alpha \in I$, $\vert\alpha\vert = r$, and if $x_i \mid x^\alpha$ then $x_{i+1}\dfrac{x^\alpha}{x_i} = \up_i(x^\alpha) \in I$. We will denote by $\{I_r\}$ the Borel set defined by $I$ and by $\{N(I)_r\}$ the set of monomials in $\poset\setminus\{I_r\}$.
\end{remark}

After having seen that a Borel ideal defines a Borel set, now we will consider the contrary perspective, deducing properties of a monomial ideal generated by a Borel set, especially concerning its Hilbert polynomial. There are some well-known facts: for example given a Borel set $B \subset \poset$ such that $\vert \poset\setminus B\vert = d \leqslant r$,we know that the ideal $I = \langle B \rangle$ has constant Hilbert polynomial $p(t) = d$ (see \cite[Theorem 3.13]{CLMR}). We introduce some definitions and notation to handle the general case.

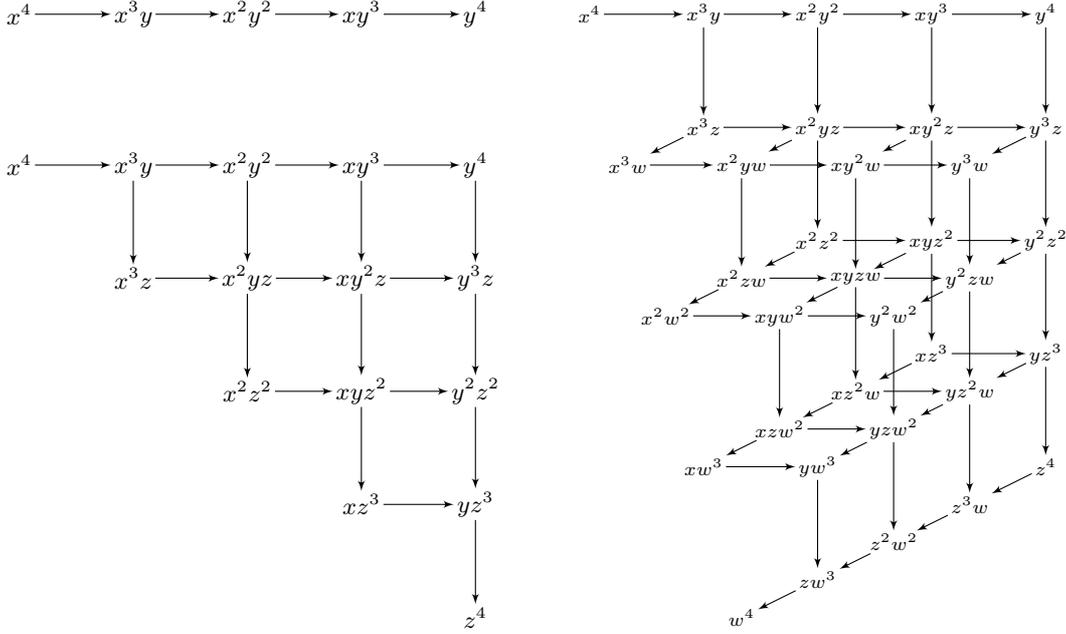
\begin{figure}[!ht]
\begin{center}
 \begin{tikzpicture}[>=latex',scale=1]
\tikzstyle{quotient}=[inner sep=1pt]
\tikzstyle{place}=[circle,draw=black,fill=black,inner sep=1pt]
\tikzstyle{place1}=[circle,draw=black,thick,inner sep=1pt]
\tikzstyle{place2}=[circle,draw=black!30,fill=black!30,inner sep=1pt]

\node (N0_1) at (-14.5,9) [quotient] {\footnotesize $x^4$};
\node (N0_2) at (-13,9) [quotient] {\footnotesize $x^3y$};
\node (N0_3) at (-11.5,9) [quotient] {\footnotesize $x^2y^2 $};
\node (N0_4) at (-10,9) [quotient] {\footnotesize $xy^3$};
\node (N0_5) at (-8.5,9) [quotient] {\footnotesize $y^4$};

\draw [->] (N0_1) -- (N0_2);
\draw [->] (N0_2) -- (N0_3);
\draw [->] (N0_3) -- (N0_4);
\draw [->] (N0_4) -- (N0_5);

\node (N1_1) at (-14.5,7) [quotient] {\footnotesize $x^4$};
\node (N1_2) at (-13,7) [quotient] {\footnotesize $x^3y$};
\node (N1_3) at (-11.5,7) [quotient] {\footnotesize $x^2y^2 $};
\node (N1_4) at (-10,7) [quotient] {\footnotesize $xy^3$};
\node (N1_5) at (-8.5,7) [quotient] {\footnotesize $y^4$};

\draw [->] (N1_1) -- (N1_2);
\draw [->] (N1_2) -- (N1_3);
\draw [->] (N1_3) -- (N1_4);
\draw [->] (N1_4) -- (N1_5);

\node (N1_6) at (-13,5.5) [quotient] {\footnotesize $x^3z$};
\node (N1_7) at (-11.5,5.5) [quotient] {\footnotesize $x^2yz $};
\node (N1_8) at (-10,5.5) [quotient] {\footnotesize $xy^2 z$};
\node (N1_9) at (-8.5,5.5) [quotient] {\footnotesize $y^3 z$};

\draw [->] (N1_2) -- (N1_6);
\draw [->] (N1_3) -- (N1_7);
\draw [->] (N1_4) -- (N1_8);
\draw [->] (N1_5) -- (N1_9);

\draw [->] (N1_6) -- (N1_7);
\draw [->] (N1_7) -- (N1_8);
\draw [->] (N1_8) -- (N1_9);

\node (N1_10) at (-11.5,4) [quotient] {\footnotesize $x^2 z^2 $};
\node (N1_11) at (-10,4) [quotient] {\footnotesize $xy z^2$};
\node (N1_12) at (-8.5,4) [quotient] {\footnotesize $y^2 z^2$};

\draw [->] (N1_7) -- (N1_10);
\draw [->] (N1_8) -- (N1_11);
\draw [->] (N1_9) -- (N1_12);

\draw [->] (N1_10) -- (N1_11);
\draw [->] (N1_11) -- (N1_12);

\node (N1_13) at (-10,2.5) [quotient] {\footnotesize $x z^3$};
\node (N1_14) at (-8.5,2.5) [quotient] {\footnotesize $y z^3$};

\draw [->] (N1_11) -- (N1_13);
\draw [->] (N1_12) -- (N1_14);

\draw [->] (N1_13) -- (N1_14);

\node (N1_15) at (-8.5,1) [quotient] {\footnotesize $z^4$};

\draw [->] (N1_14) -- (N1_15);

\node (N_1) at (-7,9) [quotient] {\tiny $x^4$};
\node (N_2) at (-5.5,9) [quotient] {\tiny $x^3y$};
\node (N_3) at (-4,9) [quotient] {\tiny $x^2y^2 $};
\node (N_4) at (-2.5,9) [quotient] {\tiny $xy^3$};
\node (N_5) at (-1,9) [quotient] {\tiny $y^4$};

\draw [->] (N_1) -- (N_2);
\draw [->] (N_2) -- (N_3);
\draw [->] (N_3) -- (N_4);
\draw [->] (N_4) -- (N_5);

\node (N_6) at (-5.5,7.5) [quotient] {\tiny $x^3z$};
\node (N_7) at (-4,7.5) [quotient] {\tiny $x^2yz $};
\node (N_8) at (-2.5,7.5) [quotient] {\tiny $xy^2 z$};
\node (N_9) at (-1,7.5) [quotient] {\tiny $y^3 z$};

\draw [->] (N_2) -- (N_6);
\draw [->] (N_3) -- (N_7);
\draw [->] (N_4) -- (N_8);
\draw [->] (N_5) -- (N_9);

\draw [->] (N_6) -- (N_7);
\draw [->] (N_7) -- (N_8);
\draw [->] (N_8) -- (N_9);

\node (N_10) at (-4,6) [quotient] {\tiny $x^2 z^2 $};
\node (N_11) at (-2.5,6) [quotient] {\tiny $xy z^2$};
\node (N_12) at (-1,6) [quotient] {\tiny $y^2 z^2$};

\draw [->] (N_7) -- (N_10);
\draw [->] (N_8) -- (N_11);
\draw [->] (N_9) -- (N_12);

\draw [->] (N_10) -- (N_11);
\draw [->] (N_11) -- (N_12);

\node (N_13) at (-2.5,4.5) [quotient] {\tiny $x z^3$};
\node (N_14) at (-1,4.5) [quotient] {\tiny $y z^3$};

\draw [->] (N_11) -- (N_13);
\draw [->] (N_12) -- (N_14);

\draw [->] (N_13) -- (N_14);

\node (N_15) at (-1,3) [quotient] {\tiny $z^4$};

\draw [->] (N_14) -- (N_15);

\node (M_6) at (-6.5,7) [quotient] {\tiny $x^3w$};
\node (M_7) at (-5,7) [quotient] {\tiny $x^2yw $};
\node (M_8) at (-3.5,7) [quotient] {\tiny $xy^2 w$};
\node (M_9) at (-2,7) [quotient] {\tiny $y^3 w$};

\draw [->] (M_6) -- (M_7);
\draw [->] (M_7) -- (M_8);
\draw [->] (M_8) -- (M_9);

\node (M_10) at (-5,5.5) [quotient] {\tiny $x^2 z w $};
\node (M_11) at (-3.5,5.5) [quotient] {\tiny $xy z w$};
\node (M_12) at (-2,5.5) [quotient] {\tiny $y^2 z w$};

\draw [->] (M_7) -- (M_10);
\draw [->] (M_8) -- (M_11);
\draw [->] (M_9) -- (M_12);

\draw [->] (M_10) -- (M_11);
\draw [->] (M_11) -- (M_12);

\node (M_13) at (-3.5,4) [quotient] {\tiny $x z^2 w$};
\node (M_14) at (-2,4) [quotient] {\tiny $y z^2 w$};

\draw [->] (M_11) -- (M_13);
\draw [->] (M_12) -- (M_14);

\draw [->] (M_13) -- (M_14);

\node (M_15) at (-2,2.5) [quotient] {\tiny $z^3 w$};

\draw [->] (M_14) -- (M_15);

\draw [->] (N_6) -- (M_6);
\draw [->] (N_7) -- (M_7);
\draw [->] (N_8) -- (M_8);
\draw [->] (N_9) -- (M_9);
\draw [->] (N_10) -- (M_10);
\draw [->] (N_11) -- (M_11);
\draw [->] (N_12) -- (M_12);
\draw [->] (N_13) -- (M_13);
\draw [->] (N_14) -- (M_14);
\draw [->] (N_15) -- (M_15);

\node (P_10) at (-6,5) [quotient] {\tiny $x^2 w^2$};
\node (P_11) at (-4.5,5) [quotient] {\tiny $xy w^2$};
\node (P_12) at (-3,5) [quotient] {\tiny $y^2 w^2$};

\draw [->] (P_10) -- (P_11);
\draw [->] (P_11) -- (P_12);

\node (P_13) at (-4.5,3.5) [quotient] {\tiny $x z w^2$};
\node (P_14) at (-3,3.5) [quotient] {\tiny $y z w^2$};

\draw [->] (P_11) -- (P_13);
\draw [->] (P_12) -- (P_14);

\draw [->] (P_13) -- (P_14);

\node (P_15) at (-3,2) [quotient] {\tiny $z^2 w^2$};

\draw [->] (P_14) -- (P_15);

\draw [->] (M_10) -- (P_10);
\draw [->] (M_11) -- (P_11);
\draw [->] (M_12) -- (P_12);
\draw [->] (M_13) -- (P_13);
\draw [->] (M_14) -- (P_14);
\draw [->] (M_15) -- (P_15);

\node (Q_13) at (-5.5,3) [quotient] {\tiny $xw^3$};
\node (Q_14) at (-4,3) [quotient] {\tiny $yw^3$};

\draw [->] (Q_13) -- (Q_14);

\node (Q_15) at (-4,1.5) [quotient] {\tiny $zw^3$};

\draw [->] (Q_14) -- (Q_15);

\draw [->] (P_13) -- (Q_13);
\draw [->] (P_14) -- (Q_14);
\draw [->] (P_15) -- (Q_15);

\node (R_15) at (-5,1) [quotient] {\tiny $w^4$};

\draw [->] (Q_15) -- (R_15);


\end{tikzpicture}
\end{center}
\caption{\label{fig:posetInclusion} The inclusion $\mathcal{P}(1,4) \subset \mathcal{P}(2,4) \subset \mathcal{P}(3,4)$ of the monomials in $k[x,y,z,w]$, $x \geq_B y \geq_B z \geq_B w$.}
\end{figure}

\begin{definition}\label{def:minMonomio}
Given a monomial $x^\alpha \in \K[x]$, we define
\[
 \min x^\alpha = \min\left\{j\text{ s.t. } x_j \mid x^\alpha\right\} \qquad\text{and}\qquad \max x^\alpha = \max\left\{j\text{ s.t. } x_j \mid x^\alpha\right\}.
\]
Moreover given a Borel set $B \subset \poset$, $N = \poset\setminus B$, for $j=0,\ldots,n-1$ let
\begin{itemize}
 \item $B_j = B \cap \mathcal{P}(n-j,r) = \{x^\alpha \in B\ \vert\ \min x^\alpha \geqslant j\}$;
 \item $N_j = N \cap \mathcal{P}(n-j,r) = \{x^\alpha \notin B\ \vert\ \min x^\alpha \geqslant j\}$.
\end{itemize}
Obviously $B_j \cup N_j = \mathcal{P}(n-j,r)$ and $B_{n-1} \subset \ldots \subset B_0$, $N_{n-1} \subset \ldots \subset N_0$ from the natural inclusion $\mathcal{P}(1,r) \subset \ldots \subset \mathcal{P}(i,r) \subset \ldots \subset \mathcal{P}(n,r)$ (see an example in Figure \ref{fig:posetInclusion}).
\end{definition}

The algorithm discussed in Section 5 of \cite{CLMR} allows to deduce the following property.
\begin{proposition}\label{prop:BorelSetHP}
Let $B$ be a Borel set in the poset $\mathcal{P}(n,r)$, $N$ defined as above and let $p(t)$ be a polynomial of degree $< n$ and such that $\vert N_j\vert = \Delta^j p(r),\ \forall\ j=0,\ldots,n-1$. Then the quotient ring $\K[x]/I$ defined by the ideal $I = \left\langle B \right\rangle$ has Hilbert polynomial equal to $p(t)$.
\begin{proof}
 The algorithm exposed in \cite{CLMR} compute all the Borel-fixed ideals in $\K[x]$ generated in degree $r$ with assigned Hilbert polynomial $p(t)$ by imposing in all the possible ways $\vert N_j\vert = \Delta^j p(r), \forall\ j=0,\ldots,n-1$.
\end{proof}
\end{proposition}

The previous proposition is useful because it gives a concrete method, shown in the following example, to recover the Hilbert polynomial of an ideal generated by a Borel set.
\begin{example}
 Let us consider the Borel set $B = \{x^4,x^3y,x^2y^2,xy^3,x^3z,x^2yz,x^2z^2\} \in \mathcal{P}(2,4)$ (Figure \ref{fig:ExampleBorel}\subref{fig:ExampleBorel_b}) and let us try to find a polynomial with the property as in Theorem \ref{prop:BorelSetHP}. First of all we easily compute $N_0 = \{y^4,xy^2z,y^3z,xyz^2,y^2z^2,xz^3,yz^3,z^4\}$ and $N_1 = \{y^4\}$. Then the degree of the polynomial $p(t)$ we are looking for has to be 1 so that $\Delta p(t)$ has to be constant and equal to $\vert N_1\vert = 1$. $p(t)$ will be of the type $t+c$ and imposing $p(4) = 4 + c = \vert N_0\vert = 8$ we find $p(t) = t + 4$ as expected.
\end{example}

We conclude this section extending some definition given in \cite{CLMR}.
\begin{definition}
 Let $B$ be a Borel set in $\poset$ and $N = \poset\setminus B$ as before.
\begin{itemize}
 \item $x^\alpha$ is a \emph{minimal element} of $B$ if $x^\alpha \in B$ and $\down_i(x^\alpha) \in N$ for all admissible $\down_i$;
 \item $x^\alpha$ is a \emph{minimal element} of $B_j$ if $x^\alpha \in B_j$ and $\down_i(x^\alpha) \in N_j$ for all admissible $\down_i,\ i > j$;
 \item $x^\beta$ is a \emph{maximal element} of $N$ if $x^\beta \in N$ and $\up_i(x^\alpha) \in B$ for all admissible $\up_i$;
 \item $x^\beta$ is a \emph{minimal element} of $N_j$ if $x^\beta \in N_j$ and $\up_i(x^\beta) \in B_j$ for all admissible $\up_i$.
\end{itemize}
\end{definition}

Keeping in mind the partial order $\leq_B$ and the inclusion $B_{j} \subset B_{j-1},\ N_j \subset N_{j-1}$, it is easy to see that a minimal element in $B_{j}$ is not necessary a minimal elements in $B_{j-1}$ because considering one more smaller variable we can obtain smaller monomials, whereas for the same reason a maximal element of $N_j$ remains maximal also in $N_{j-1}$. For example, let us consider the Borel set $B = \{x^3,x^2y,x^2z\} \in \mathcal{P}(2,3)$: $x^2y$ is minimal in $B_1$ because $\down_x(x^2 y) = xy^2 \in N_1$, but it is not minimal in $B_0$, because $\down_{y}(x^2y) = x^2 z \in B$. On the other hand $xy^2$ is maximal both in $N_1$ and in $N_0$.

\begin{remark}\label{rk:AddRemove}
 Let $B \subset \poset$ be a Borel set and let $x^\alpha$ and $x^\beta$ be a minimal and a maximal monomial. By definition both $B\setminus\{x^\alpha\}$ and $B\cup\{x^\beta\}$ are still Borel sets. Moreover $x^\alpha$ will be a maximal element of $\poset \setminus \left(B \setminus \{x^\alpha\}\right)$ and $x^\beta$ will be a minimal element of $B \cup \{x^\beta\}$.
\end{remark}

\begin{remark}\label{rk:minimalMinimum}
 Let $\preceq$ be a term ordering such that $x_n \succ \cdots \succ x_0$ and let $B$ be a Borel set in $\poset$. Since $\preceq$ is a refinement of the partial order $\leq_B$, the minimum monomial $x^\alpha = \min_\preceq \{x^\gamma \in B_j\}$ is a minimal element of $B_j$ and the maximum monomial $x^\beta = \max_\preceq \{x^\delta \in N_j\}$ is a maximal element of $N_j$.

On the contrary, if $\{x^{\alpha_1},\ldots,x^{\alpha_t}\}$ is the set of the minimal elements in $B_j$, the minimum monomial $x^\alpha \in B_j$ w.r.t. any term order $\preceq$ (refinement of $\leq_B$) will belong to this set and analogously the set $\{x^{\beta_1},\ldots,x^{\beta_s}\}$ of the maximal elements in $N_j$ will contain the maximum $\max_{\preceq}\{x^\delta \in N_j\}$.
\end{remark}

\section{The zero-dimensional case}\label{sec:0case}
Let us consider a Borel set $B \subset \poset$ such that $\vert N_0\vert = d \leqslant r$ and $N_j = \emptyset, \forall\ j > 0$. By Proposition \ref{prop:BorelSetHP} we now that the ideal $I = \langle B\rangle$ defines a subscheme of $\PP^n$ with constant Hilbert polynomial $p(t) = d$. Moreover let us suppose that there exist two monomials $x^\alpha,x^\beta \in \poset$, such that $x^\alpha$ is a minimal element in $B$ such that $\min x^\alpha = 0$ and $x^\beta$ is a maximal one in $N$ (surely $\min x^\beta = 0$). If the set $\widetilde{B} = B  \cup \{x^\beta\}\setminus\{x^\alpha\} $ is still Borel, since $\vert\widetilde{N}_0\vert = \vert N_0\vert$, the ideal $\widetilde{I} = \langle\widetilde{B}\rangle$ has the same Hilbert polynomial of $I$.

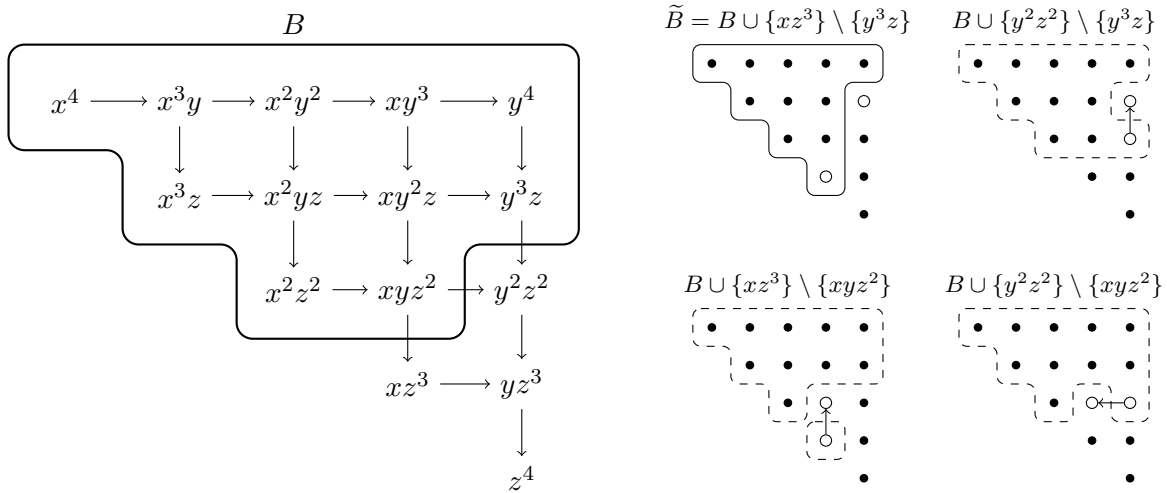
\begin{figure}[!ht] 
\begin{center}
 \begin{tikzpicture}
\tikzstyle{quotient}=[]
\node (M_1) at (0.5,5) [quotient] {$x^4$};
\node (M_2) at (2,5) [quotient] {$x^3y$};
\node (M_3) at (3.5,5) [quotient] {$x^2y^2$};
\node (M_4) at (5,5) [quotient] {$xy^3$};
\node (M_5) at (6.5,5) [quotient] {$y^4$};

\node (M_6) at (2,3.75) [quotient] {$x^3z$};
\node (M_7) at (3.5,3.75) [quotient] {$x^2yz$};
\node (M_8) at (5,3.75) [quotient] {$xy^2z$};
\node (M_9) at (6.5,3.75) [quotient] {$y^3z$};

\node (M_10) at (3.5,2.5) [quotient] {$x^2z^2$};
\node (M_11) at (5,2.5) [quotient] {$xyz^2$};
\node (M_12) at (6.5,2.5) [quotient] {$y^2z^2$};

\node (M_13) at (5,1.25) [quotient] {$xz^3$};
\node (M_14) at (6.5,1.25) [quotient] {$yz^3$};

\node (M_15) at (6.5,0) [quotient] {$z^4$};

\draw [->] (M_1) -- (M_2);
\draw [->] (M_2) -- (M_3);
\draw [->] (M_3) -- (M_4);
\draw [->] (M_4) -- (M_5);

\draw [->] (M_2) -- (M_6);
\draw [->] (M_3) -- (M_7);
\draw [->] (M_4) -- (M_8);
\draw [->] (M_5) -- (M_9);

\draw [->] (M_6) -- (M_7);
\draw [->] (M_7) -- (M_8);
\draw [->] (M_8) -- (M_9);

\draw [->] (M_7) -- (M_10);
\draw [->] (M_8) -- (M_11);
\draw [->] (M_9) -- (M_12);

\draw [->] (M_10) -- (M_11);
\draw [->] (M_11) -- (M_12);

\draw [->] (M_11) -- (M_13);
\draw [->] (M_12) -- (M_14);

\draw [->] (M_13) -- (M_14);
\draw [->] (M_14) -- (M_15);

\draw [rounded corners=6pt,thick] (4.5,1.85) -- (5.75,1.85) -- (5.75,3.1) -- (7.25,3.1) -- (7.25,5.75) --  node[above] {$B$} (-0.25,5.75) -- (-0.25,4.35) -- (1.25,4.35) -- (1.25,3.1) -- (2.75,3.1) -- (2.75,1.85) -- cycle;

\tikzstyle{place}=[circle,draw=black,fill=black,inner sep=1pt]
\tikzstyle{place1}=[circle,draw=black,inner sep=1.5pt]

\node at (9,2) [place] {};
\node at (9.5,2) [place] {};
\node at (10,2) [place] {};
\node at (10.5,2) [place] {};
\node at (11,2) [place] {};

\node at (9.5,1.5) [place] {};
\node at (10,1.5) [place] {};
\node at (10.5,1.5) [place] {};
\node at (11,1.5) [place] {};

\node at (10,1) [place] {};
\node (W_3) at (10.5,1) [place1] {};
\node at (11,1) [place] {};

\node (W_4) at (10.5,0.5) [place1] {};
\node at (11,0.5) [place] {};

\draw [->] (W_4) -- (W_3);

\node at (11,0) [place] {};

\draw[rounded corners=3pt,dashed] (8.75,2.25) -- node[above]{\footnotesize $B \cup\{xz^3\}\setminus\{xyz^2\}$} (11.25,2.25) -- (11.25,1.25) -- (10.25,1.25) -- (10.25,0.75) -- (9.75,0.75) -- (9.75,1.25) -- (9.25,1.25) -- (9.25,1.75) -- (8.75,1.75) -- cycle;

\draw[rounded corners=3pt,dashed] (10.25,0.75) -- (10.75,0.75) -- (10.75,0.25) -- (10.25,0.25) -- cycle;

\node at (12.5,2) [place] {};
\node at (13,2) [place] {};
\node at (13.5,2) [place] {};
\node at (14,2) [place] {};
\node at (14.5,2) [place] {};

\node at (13,1.5) [place] {};
\node at (13.5,1.5) [place] {};
\node at (14,1.5) [place] {};
\node at (14.5,1.5) [place] {};

\node at (13.5,1) [place] {};
\node (W_5) at (14,1) [place1] {};
\node (W_6) at (14.5,1) [place1] {};

\draw [->] (W_6) -- (W_5);

\node at (14,0.5) [place] {};
\node at (14.5,0.5) [place] {};

\node at (14.5,0) [place] {};

\draw[rounded corners=3pt,dashed] (12.25,2.25) -- node[above]{\footnotesize $B \cup\{y^2z^2\}\setminus\{xyz^2\}$} (14.75,2.25) -- (14.75,0.75) -- (14.25,0.75) -- (14.25,1.25) -- (13.75,1.25) -- (13.75,0.75) -- (13.25,0.75) -- (13.25,1.25) -- (12.75,1.25) -- (12.75,1.75) -- (12.25,1.75) -- cycle;

\node at (9,5.5) [place] {};
\node at (9.5,5.5) [place] {};
\node at (10,5.5) [place] {};
\node at (10.5,5.5) [place] {};
\node at (11,5.5) [place] {};

\node at (9.5,5) [place] {};
\node at (10,5) [place] {};
\node at (10.5,5) [place] {};
\node at (11,5) [place1] {};

\node at (10,4.5) [place] {};
\node at (10.5,4.5) [place] {};
\node at (11,4.5) [place] {};

\node at (10.5,4) [place1] {};
\node at (11,4) [place] {};

\node at (11,3.5) [place] {};

\draw[rounded corners=3pt] (8.75,5.75) -- node[above]{\footnotesize $\widetilde{B}= B \cup\{xz^3\}\setminus\{y^3z\}$} (11.25,5.75) -- (11.25,5.25) -- (10.75,5.25) -- (10.75,3.75) -- (10.25,3.75) -- (10.25,4.25) -- (9.75,4.25) -- (9.75,4.75) -- (9.25,4.75) -- (9.25,5.25) -- (8.75,5.25) -- cycle;

\node at (12.5,5.5) [place] {};
\node at (13,5.5) [place] {};
\node at (13.5,5.5) [place] {};
\node at (14,5.5) [place] {};
\node at (14.5,5.5) [place] {};

\node at (13,5) [place] {};
\node at (13.5,5) [place] {};
\node at (14,5) [place] {};
\node(W_1) at (14.5,5) [place1] {};

\node at (13.5,4.5) [place] {};
\node at (14,4.5) [place] {};
\node(W_2) at (14.5,4.5) [place1] {};

\draw[->] (W_2) -- (W_1);
\node at (14,4) [place] {};
\node at (14.5,4) [place] {};

\node at (14.5,3.5) [place] {};

\draw[rounded corners=3pt,dashed] (12.25,5.75) -- node[above]{\footnotesize $B \cup\{y^2z^2\}\setminus\{y^3z\}$} (14.75,5.75) -- (14.75,5.25) -- (14.25,5.25) -- (14.25,4.75) -- (14.75,4.75) -- (14.75,4.25) -- (13.25,4.25) -- (13.25,4.75) -- (12.75,4.75) -- (12.75,5.25) -- (12.25,5.25) -- cycle;
\end{tikzpicture}
\end{center}
\caption{\label{fig:ExampleSwap}On the left there is the Borel set $B$ described in Example \ref{ex:ExampleSwap}. On the right there are the sets that can be obtained exchanging a minimal and a maximal monomial (the circled vertices). $\widetilde{B}$ is the only one preserving the Borel condition. In the other three cases we highlight an increasing elementary move that makes that condition false.}
\end{figure}

\begin{example}\label{ex:ExampleSwap}
Let us consider the Borel set $B = \{x^4,x^3y,x^2y^2,xy^3,y^4,x^3z,x^2yz,xy^2z,y^3z,$ $x^2z^2,xyz^2\} \subseteq \mathcal{P}(2,4)$. The monomials $xyz^2$ and $y^3z$ are minimal in $B$ and $xz^3$ and $y^2z^2$ are maximal in $N$. The only exchange that preserves the Borel condition is $\widetilde{B} = B \cup\{xz^3\}\setminus\{y^3z\}$, whereas the sets $B\cup\{y^2z^2\} \setminus\{y^3z\}$, $B \cup\{xz^3\}\setminus\{xyz^2\}$ and $B \cup\{y^2z^2\}\setminus\{xyz^2\}$ are not Borel (see Figure \ref{fig:ExampleSwap}).
\end{example}

By Remark \ref{rk:AddRemove} we know that we can add a maximal element or remove a minimal one to a Borel set preserving this property. In Example \ref{ex:ExampleSwap} we obtain sets of monomials not Borel because after adding the maximal element $x^\beta$ to $B$, $x^\alpha$ is no more minimal in $B\cup\{x^\beta\}$. This happens whenever an elementary move $\down_i$ such that $\down_i(x^\alpha) = x^\beta$ exists. 

We are ready to state the first result of this paper.

\begin{theorem}\label{th:flatZeroDim}
Let $B\subset\mathcal{P}(n,d)$ be a Borel set and $N$ be defined as before such that $\vert N_0 \vert = d$ and $N_j = \emptyset,\ \forall\ j > 0$. Moreover let us suppose that there exist a minimal monomial $x^\alpha \in B$ and a maximal $x^\beta \in N$, such that $\min x^\alpha = 0$ and $\down_j(x^\alpha) \neq x^\beta$ for every admissible $\down_j$. The ideal 
\[
 \mathcal{I} = \left\langle B\setminus\left\{x^\alpha\right\}\cup\left\{y_0 \, x^\alpha + y_1 \, x^\beta\right\}\right\rangle
\]
has Hilbert polynomial $\binom{n+t}{n}-d$, for all $[y_0:y_1] \in \PP^1$.
\begin{proof}
First of all let us remark that the ideals $\mathcal{I}\vert_{[1:0]} = \langle B\rangle$ and $\mathcal{I}\vert_{[0:1]} = \big\langle B\setminus\{x^\alpha\}\cup\{x^\beta\}\big\rangle$ have Hilbert polynomial $\binom{n+t}{n}-d$. In fact both $B$ and $\widetilde{B} = B\setminus\{x^\alpha\}\cup\{x^\beta\}$ are Borel sets and they are in the hypothesis of Proposition \ref{prop:BorelSetHP} with polynomial $p(t) = d$.

Since $\dim_k \mathcal{I}_d = \binom{n+d}{n}-d,\ \forall\ [y_0:y_1] \in \PP^1$ by construction, to obtain the thesis it suffices to prove $\dim_k \mathcal{I}_{d+1} = \binom{n+d+1}{n}-d$ and then to apply Theorem \ref{th:Gotzmann}.

Let us suppose $y_0 \neq 0$ and let us denote by $J$ the ideal $\mathcal{I}\vert_{[1:y_1]}$ and by $I = \mathcal{I}\vert_{[1:0]}$. We want to show that $ \dim_k J_{d+1} = \dim_k I_{d+1}$.
Let $\mathcal{G}$ be the standard monomial basis of $I_{d+1}$. Among the monomials of $\mathcal{G}$ there is $x^\alpha x_0$ and by the minimality of $x^\alpha$, there is no other monomial $x^\gamma \in B$ such that $x^\alpha x_0 = x^\gamma x_i$, for any $i > 0$, because this would imply $x^\gamma \leq_B x^\alpha$, since $x_i \geq_B x_0$. On the contrary for any $i > 0$, $x^\alpha x_i$ could be generated from another monomial $x^\delta \in B\setminus\{x^\alpha\}$, in fact
\begin{equation}\label{eq:otherGenerator}
 x_i\,x^\alpha = \up_{i-1}\cdot\ldots\cdot\up_1\cdot\up_0(x^\alpha)\, x_0 = x^\delta\, x_0 \qquad(\min x^\alpha = 0).
\end{equation}

We claim that $\mathcal{G}\setminus\{x^\alpha x_0\}\cup\{x^\alpha x_0 + y_1 \, x^\beta x_0\}$ is a basis of $J_{r+1}$. To prove this claim, it is sufficient to verify that for any $i > 0$, both monomials in $x_i(x^\alpha + y_1\, x^\beta)$ can be obtained from monomials in $B\setminus\{x^\alpha\}$:
\begin{equation}
x_i(x^\alpha + y_1\,x^\beta) = x^\delta\,x_0 + y_1\, \up_{i-1}\cdot\ldots\cdot\up_1\cdot\up_0(x^\beta)\, x_0 \qquad(\min x^\beta = 0),
\end{equation}
where $x^\delta$ and $\up_{i-1}\cdot\ldots\cdot\up_1\cdot\up_0(x^\beta)$ belong to $B\setminus\{x^\alpha\}$ by \eqref{eq:otherGenerator} and by the hypothesis of maximality of $x^\beta$ in $N$.

Obviously, the same reasoning holds in the case $y_1 \neq 0$ exchanging the role of $x^\alpha$ and $x^\beta$ (Remark \ref{rk:AddRemove}), so in conclusion $\forall\ [y_0 : y_1] \in \PP^1$
\[
\dim_k \mathcal{I}_{d+1} = \dim_k \left(\mathcal{I}\vert_{(1,0)}\right)_{d+1} = \dim_k \left(\mathcal{I}\vert_{(0,1)}\right)_{d+1} = \binom{n+d+1}{n}-d
\]
and by Gotzmann Persistence Theorem $\dim_k \mathcal{I}_{t} = \binom{n+t}{n}-d,\ \forall\ t \geqslant d$.
\end{proof}
\end{theorem}

\begin{example}
The ideals $I = (x^2,xy,y^3)_{\geqslant 4}$ and $\widetilde{I} = (x,y^4)_{\geqslant 4}$ defined in Example \ref{ex:ExampleSwap} are in the hypothesis of Theorem \ref{th:flatZeroDim}, so all ideals of the family 
\[
 \mathcal{I} = \left\langle\left\{I_4\right\}\setminus\left\{y^3 z\right\}\cup\left\{u_0\, y^3 z + u_1\, xz^3\right\}\right\rangle,\quad [u_0:u_1]\in\PP^1,
\]
have Hilbert polynomial $q(t) = \binom{t+2}{2} - 4$. 
\end{example}

\section{The general case}\label{sec:genCase}
We would like to extend the technique to Borel ideals with Hilbert polynomial $p(t)$ of any degree. The point is to understand how to substitute monomials in a Borel set $B$ such that $N_j \neq \emptyset$, for $j > 0$, in order to obtain another Borel set $\widetilde{B}$, with the same Hilbert polynomial, that is $\vert B_j\vert = \vert\widetilde{B}_j\vert,\ \forall\ j$. 
The first idea is to exchange a minimal monomial $x^\alpha \in B$ with a maximal $x^\beta \in N$, but whenever $\min x^\alpha \neq \min x^\beta$ this exchange will not preserve the Hilbert polynomial.

Let us start introducing some further notation.
\begin{notation}
Given any Borel set $B$, we denote by $\mathfrak{m}(B_j)$ the set of all the minimal elements of $B_j$ contained in $B_j \setminus B_{j+1}$ and by $\mathfrak{M}(N_j)$ the set of all the maximal elements of $N_j$ contained in $N_j \setminus N_{j+1}$. When the Borel set is defined by an ideal, we will write $\mathfrak{m}\{I_r\}_j$ and $\mathfrak{M}\{N(I)_r\}_j$ instead of $\mathfrak{m}\big(\{I_r\}_j\big)$ and $\mathfrak{M}\big(\{N(I)_r\}_j\big)$.
\end{notation}

Therefore a second idea could be to exchange a monomial $x^\alpha$ in $\mathfrak{m}(B_i)$ and a monomial $x^\beta$ in $\mathfrak{M}(N_i)$. In this way, the cardinality of the sets $B_i$ and $\widetilde{B}_i$ is preserved, but whenever $\down_{i}(x^\alpha)$ belongs to $B$, swapping $x^\alpha$ and $x^\beta$ we do not obtain a Borel set. This fact suggests that the general case is more complicated and that we have to swap more monomials than the minimal one and the maximal one considered.

\begin{definition}\label{def:decreasingMoves}
 Let $B \subset \poset$ be a Borel set, $N = \poset\setminus B$ and let $x^\alpha$ be a monomial in $\mathfrak{m}(B_j)$ and $x^\beta$ a monomial in $\mathfrak{M}(N_j)$. We define the family of the admissible \emph{decreasing moves} w.r.t. $x^\alpha$ as the set
 \begin{equation}
  \mathcal{F}_\alpha = \left\{ F = (\down_{l_s})^{\lambda_s} \cdot \ldots \cdot (\down_{l_1})^{\lambda_1},\ j \geqslant l_1 > \ldots > l_s\  \big\vert\ F(x^\alpha) \in B \right\} \cup \{\mathrm{id}\},
 \end{equation}
where $\mathrm{id}$ is the \lq\lq identity move\rq\rq\ that leaves fixed each monomial.
We will say that $\mathcal{F}_\alpha$ is \emph{Borel-consistent} w.r.t. $x^\beta$ if every $F \in \mathcal{F}_\alpha$ is admissible also for $x^\beta$ and if $\up_i\big(F(x^\beta)\big) \in B$ for all admissible elementary move $\up_i,\ i \geqslant j$. In this case we define the two sets of monomials
\[
 \mathcal{F}_\alpha(x^\alpha) = \left\{ F(x^\alpha) \ \vert\ F \in \mathcal{F}_\alpha \right\}\subset B\qquad\text{and}\qquad
 \mathcal{F}_{\alpha}(x^\beta) = \left\{ F(x^\beta) \ \vert\ F \in \mathcal{F}_\alpha \right\}\subset N.
\]
\end{definition}
 
\begin{lemma}\label{lem:PreservingBorel}
 Let $B,N,x^\alpha,x^\beta$ be defined as in Definition \ref{def:decreasingMoves}. If for any admissible $\down_k$, $k > j$, $\down_{k}(x^\alpha) \neq x^\beta$ and if the family of decreasing moves $\mathcal{F}_\alpha$ is Borel-consistent w.r.t. $x^\beta$, then $\widetilde{B} = B \cup \mathcal{F}_\alpha(x^\beta) \setminus \mathcal{F}_\alpha(x^\alpha)$ is a Borel set and $\vert B_i \vert = \vert \widetilde{B}_i\vert,\ \forall\ i$.
\begin{proof}
The monomials in $\mathcal{F}_\alpha(x^\alpha)$ are in the \lq\lq internal border\rq\rq\ of $B$ in the sense that for any $\down_i,\ i > j$ and for any $F \in \mathcal{F}_\alpha$, $\down_i\big(F(x^\alpha)\big)$ does not belong to $B$, because $\down_i\big(F(x^\alpha)\big) = F \big( \down_i(x^\alpha)\big)$ and $F\big(\down_i (x^\alpha)\big) \leq_B  \down_i (x^\alpha) \notin B$. On the other hand the monomials in $\mathcal{F}_\alpha(x^\beta)$ are in the \lq\lq external border\rq\rq\ of $B$ by definition, so swapping these two sets of monomials preserves the Borel condition.
 Finally since for every $F \in \mathcal{F}_\alpha$, $\min F(x^\alpha) = \min F(x^\beta)$ also the condition $\vert B_i \vert = \vert \widetilde{B}_i\vert,\ \forall\ i$ holds.
\end{proof}
\end{lemma}

\begin{example}\label{ex:breakPreserveBorel}
Let us consider the poset $\mathcal{P}(3,4)$ (represented in Figure \ref{fig:posetInclusion}).
\begin{enumerate}
 \item Let $B$ be the Borel set defined by $(x,y^4)_4$ (Figure \ref{fig:ExampleSubstMonom}\subref{fig:ExampleSubstMonom_1}). $xz^3$ is a minimal monomial in $B_1$ and $y^3z$ a maximal one in $N_1$. The set of decreasing moves starting from $xz^3$ is $\mathcal{F} = \{\mathrm{id},\down_z,(\down_z)^2,(\down_z)^3\}$, but only the first two are admissible w.r.t. $y^3z$, so we cannot swap these two monomials.
 \item Let now consider the Borel set $\widetilde{B}$ defined by $(x^2,xy^2,y^3,xyz^2)_4$ (Figure \ref{fig:ExampleSubstMonom}\subref{fig:ExampleSubstMonom_2}). $y^3z$ is a minimal monomial in $\widetilde{B}_1$ and $xz^3$ a maximal one in $\widetilde{N}_1$. The set of decreasing moves starting from $y^3z$ is $\widetilde{\mathcal{F}} = \{\mathrm{id},\down_z\}$, also admissible w.r.t. $xz^3$, but $\up_z \big(\down_z(xz^3)\big) = xyzw \notin \widetilde{B}$. Exchanging $\widetilde{\mathcal{F}}(y^3z)$ with $\widetilde{\mathcal{F}}(xz^3)$ would break the Borel condition.
 \item Finally let $\widehat{B}$ be the Borel set defined by $(x^2,xy,xz^2,y^4)_4$ (Figure \ref{fig:ExampleSubstMonom}\subref{fig:ExampleSubstMonom_3}). $xz^3$ is minimal in $\widehat{B}_1$ and $y^3z$ maximal in $\widehat{N}_1$. $\widehat{\mathcal{F}} = \{\mathrm{id},\down_z\}$, $\widehat{\mathcal{F}}(xz^3)$ and $\widehat{\mathcal{F}}(y^3z)$ respect the hypothesis of Lemma \ref{lem:PreservingBorel}, so swapping the two sets of monomials we obtain another Borel set, indeed defined by $(x^2,xy,y^3)_4$.
\end{enumerate}
\begin{figure}[!ht] 
\begin{center}\captionsetup[subfloat]{width=4.5cm}
\subfloat[][$B = \{(x,y^4)_4\}$]{\label{fig:ExampleSubstMonom_1}
\begin{tikzpicture}[scale=0.65]
\tikzstyle{quotient}=[]
\tikzstyle{place}=[circle,draw=black,fill=black,inner sep=1pt]
\tikzstyle{place1}=[circle,draw=black,thick,inner sep=1pt]
\tikzstyle{place2}=[circle,draw=black!10,fill=black!30,inner sep=1pt]

\node (N1_1) at (-15,0) [place] {};
\node (N1_2) at (-13.5,0) [place] {};
\node (N1_3) at (-12,0) [place] {};
\node (N1_4) at (-10.5,0) [place] {};
\node (N1_5) at (-9,0) [place] {};

\draw [thin,color=black!25] (N1_1) -- (N1_2);
\draw [thin,color=black!25] (N1_2) -- (N1_3);
\draw [thin,color=black!25] (N1_3) -- (N1_4);
\draw [thin,color=black!25] (N1_4) -- (N1_5);

\node (N1_6) at (-13.5,-1.5) [place] {};
\node (N1_7) at (-12,-1.5) [place] {};
\node (N1_8) at (-10.5,-1.5) [place] {};
\node (N1_9) at (-9,-1.5) [place1] {};

\draw [thin,color=black!25] (N1_2) -- (N1_6);
\draw [thin,color=black!25] (N1_3) -- (N1_7);
\draw [thin,color=black!25] (N1_4) -- (N1_8);
\draw [thin,color=black!25] (N1_5) -- (N1_9);

\draw [thin,color=black!25] (N1_6) -- (N1_7);
\draw [thin,color=black!25] (N1_7) -- (N1_8);
\draw [thin,color=black!25] (N1_8) -- (N1_9);

\node (N1_10) at (-12,-3) [place2] {};
\node (N1_11) at (-10.5,-3) [place] {};
\node (N1_12) at (-9,-3) [place1] {};

\draw [thin,color=black!25] (N1_7) -- (N1_10);
\draw [thin,color=black!25] (N1_8) -- (N1_11);
\draw [thin,color=black!25] (N1_9) -- (N1_12);

\draw [thin,color=black!25] (N1_10) -- (N1_11);
\draw [thin,color=black!25] (N1_11) -- (N1_12);

\node (N1_13) at (-10.5,-4.5) [place] {};
\node (N1_14) at (-9,-4.5) [place1] {};

\draw [thin,color=black!25] (N1_11) -- (N1_13);
\draw [thin,color=black!25] (N1_12) -- (N1_14);

\draw [thin,color=black!25] (N1_13) -- (N1_14);

\node (N1_15) at (-9,-6) [place1] {};

\draw [thin,color=black!25] (N1_14) -- (N1_15);

\node (M1_6) at (-14.5,-2) [place] {};
\node (M1_7) at (-13,-2) [place] {};
\node (M1_8) at (-11.5,-2) [place] {};
\node (M1_9) at (-10,-2) [place1] {};

\draw [thin,color=black!25] (M1_6) -- (M1_7);
\draw [thin,color=black!25] (M1_7) -- (M1_8);
\draw [thin,color=black!25] (M1_8) -- (M1_9);

\node (M1_10) at (-13,-3.5) [place] {};
\node (M1_11) at (-11.5,-3.5) [place] {};
\node (M1_12) at (-10,-3.5) [place1] {};

\draw [thin,color=black!25] (M1_7) -- (M1_10);
\draw [thin,color=black!25] (M1_8) -- (M1_11);
\draw [thin,color=black!25] (M1_9) -- (M1_12);

\draw [thin,color=black!25] (M1_10) -- (M1_11);
\draw [thin,color=black!25] (M1_11) -- (M1_12);

\node (M1_13) at (-11.5,-5) [place] {};
\node (M1_14) at (-10,-5) [place1] {};

\draw [thin,color=black!25] (M1_11) -- (M1_13);
\draw [thin,color=black!25] (M1_12) -- (M1_14);

\draw [thin,color=black!25] (M1_13) -- (M1_14);

\node (M1_15) at (-10,-6.5) [place1] {};

\draw [thin,color=black!25] (M1_14) -- (M1_15);

\draw [thin,color=black!25] (N1_6) -- (M1_6);
\draw [thin,color=black!25] (N1_7) -- (M1_7);
\draw [thin,color=black!25] (N1_8) -- (M1_8);
\draw [thin,color=black!25] (N1_9) -- (M1_9);
\draw [thin,color=black!25] (N1_10) -- (M1_10);
\draw [thin,color=black!25] (N1_11) -- (M1_11);
\draw [thin,color=black!25] (N1_12) -- (M1_12);
\draw [thin,color=black!25] (N1_13) -- (M1_13);
\draw [thin,color=black!25] (N1_14) -- (M1_14);
\draw [thin,color=black!25] (N1_15) -- (M1_15);

\node (P1_10) at (-14,-4) [place] {};
\node (P1_11) at (-12.5,-4) [place] {};
\node (P1_12) at (-11,-4) [place1] {};

\draw [thin,color=black!25] (P1_10) -- (P1_11);
\draw [thin,color=black!25] (P1_11) -- (P1_12);

\node (P1_13) at (-12.5,-5.5) [place] {};
\node (P1_14) at (-11,-5.5) [place1] {};

\draw [thin,color=black!25] (P1_11) -- (P1_13);
\draw [thin,color=black!25] (P1_12) -- (P1_14);

\draw [thin,color=black!25] (P1_13) -- (P1_14);

\node (P1_15) at (-11,-7) [place1] {};

\draw [thin,color=black!25] (P1_14) -- (P1_15);

\draw [thin,color=black!25] (M1_10) -- (P1_10);
\draw [thin,color=black!25] (M1_11) -- (P1_11);
\draw [thin,color=black!25] (M1_12) -- (P1_12);
\draw [thin,color=black!25] (M1_13) -- (P1_13);
\draw [thin,color=black!25] (M1_14) -- (P1_14);
\draw [thin,color=black!25] (M1_15) -- (P1_15);

\node (Q1_13) at (-13.5,-6) [place] {};
\node (Q1_14) at (-12,-6) [place1] {};

\draw [thin,color=black!25] (Q1_13) -- (Q1_14);

\node (Q1_15) at (-12,-7.5) [place1] {};

\draw [thin,color=black!25] (Q1_14) -- (Q1_15);

\draw [thin,color=black!25] (P1_13) -- (Q1_13);
\draw [thin,color=black!25] (P1_14) -- (Q1_14);
\draw [thin,color=black!25] (P1_15) -- (Q1_15);

\node (R1_15) at (-13,-8) [place1] {};

\draw [thin,color=black!25] (Q1_15) -- (R1_15);

\draw [dashed] (N1_9) circle (0.4);
\draw [dashed] (M1_9) circle (0.4);
\node (P1_9) at (-11,-2.5) [quotient] {};
\node (Q1_9) at (-12,-3) [quotient] {};
\draw [dashed] (P1_9) circle (0.4);
\draw [dashed] (Q1_9) circle (0.4);
\draw [dashed] (N1_13) circle (0.4);
\draw [dashed] (M1_13) circle (0.4);
\draw [dashed] (P1_13) circle (0.4);
\draw [dashed] (Q1_13) circle (0.4);
\draw [<->] (P1_13) to [bend left=20] node[fill=white,inner sep=1pt]{\textbf{?}} (P1_9);
\draw [<->] (Q1_13) to [bend left=20] node[fill=white,inner sep=1pt]{\textbf{?}} (Q1_9);

\draw [dotted] (-8.6,-1.5) -- (-11.6,-3);
\draw [dotted] (-9,-1.9) -- (-12,-3.4);
\draw [dotted] (-9,-1.1) -- (-12,-2.6);

\draw [dotted] (-10.1,-4.5) -- (-13.1,-6);
\draw [dotted] (-10.5,-4.1) -- (-13.5,-5.6);
\draw [dotted] (-10.5,-4.9) -- (-13.5,-6.4);
\end{tikzpicture}
}\quad
\subfloat[][$\widetilde{B}=\{(x^2,xy^2,y^3,xyz^2)_4\}$]{\label{fig:ExampleSubstMonom_2}
\begin{tikzpicture}[scale=0.65]
\tikzstyle{quotient}=[]
\tikzstyle{place}=[circle,draw=black,fill=black,inner sep=1pt]
\tikzstyle{place1}=[circle,draw=black,thick,inner sep=1pt]
\tikzstyle{place2}=[circle,draw=black!10,fill=black!30,inner sep=1pt]

\node (N3_1) at (-7,0) [place] {};
\node (N3_2) at (-5.5,0) [place] {};
\node (N3_3) at (-4,0) [place] {};
\node (N3_4) at (-2.5,0) [place] {};
\node (N3_5) at (-1,0) [place] {};

\draw [thin,color=black!25] (N3_1) -- (N3_2);
\draw [thin,color=black!25] (N3_2) -- (N3_3);
\draw [thin,color=black!25] (N3_3) -- (N3_4);
\draw [thin,color=black!25] (N3_4) -- (N3_5);

\node (N3_6) at (-5.5,-1.5) [place] {};
\node (N3_7) at (-4,-1.5) [place] {};
\node (N3_8) at (-2.5,-1.5) [place] {};
\node (N3_9) at (-1,-1.5) [place] {};

\draw [thin,color=black!25] (N3_2) -- (N3_6);
\draw [thin,color=black!25] (N3_3) -- (N3_7);
\draw [thin,color=black!25] (N3_4) -- (N3_8);
\draw [thin,color=black!25] (N3_5) -- (N3_9);

\draw [thin,color=black!25] (N3_6) -- (N3_7);
\draw [thin,color=black!25] (N3_7) -- (N3_8);
\draw [thin,color=black!25] (N3_8) -- (N3_9);

\node (N3_10) at (-4,-3) [place] {};
\node (N3_11) at (-2.5,-3) [place] {};
\node (N3_12) at (-1,-3) [place1] {};

\draw [thin,color=black!25] (N3_7) -- (N3_10);
\draw [thin,color=black!25] (N3_8) -- (N3_11);
\draw [thin,color=black!25] (N3_9) -- (N3_12);

\draw [thin,color=black!25] (N3_10) -- (N3_11);
\draw [thin,color=black!25] (N3_11) -- (N3_12);

\node (N3_13) at (-2.5,-4.5) [place1] {};
\node (N3_14) at (-1,-4.5) [place1] {};

\draw [thin,color=black!25] (N3_11) -- (N3_13);
\draw [thin,color=black!25] (N3_12) -- (N3_14);

\draw [thin,color=black!25] (N3_13) -- (N3_14);

\node (N3_15) at (-1,-6) [place1] {};

\draw [thin,color=black!25] (N3_14) -- (N3_15);

\node (M3_6) at (-6.5,-2) [place] {};
\node (M3_7) at (-5,-2) [place] {};
\node (M3_8) at (-3.5,-2) [place] {};
\node (M3_9) at (-2,-2) [place] {};

\draw [thin,color=black!25] (M3_6) -- (M3_7);
\draw [thin,color=black!25] (M3_7) -- (M3_8);
\draw [thin,color=black!25] (M3_8) -- (M3_9);

\node (M3_10) at (-5,-3.5) [place] {};
\node (M3_11) at (-3.5,-3.5) [place1] {};
\node (M3_12) at (-2,-3.5) [place1] {};

\draw [thin,color=black!25] (M3_7) -- (M3_10);
\draw [thin,color=black!25] (M3_8) -- (M3_11);
\draw [thin,color=black!25] (M3_9) -- (M3_12);

\draw [thin,color=black!25] (M3_10) -- (M3_11);
\draw [thin,color=black!25] (M3_11) -- (M3_12);

\node (M3_13) at (-3.5,-5) [place1] {};
\node (M3_14) at (-2,-5) [place1] {};

\draw [thin,color=black!25] (M3_11) -- (M3_13);
\draw [thin,color=black!25] (M3_12) -- (M3_14);

\draw [thin,color=black!25] (M3_13) -- (M3_14);

\node (M3_15) at (-2,-6.5) [place1] {};

\draw [thin,color=black!25] (M3_14) -- (M3_15);

\draw [thin,color=black!25] (N3_6) -- (M3_6);
\draw [thin,color=black!25] (N3_7) -- (M3_7);
\draw [thin,color=black!25] (N3_8) -- (M3_8);
\draw [thin,color=black!25] (N3_9) -- (M3_9);
\draw [thin,color=black!25] (N3_10) -- (M3_10);
\draw [thin,color=black!25] (N3_11) -- (M3_11);
\draw [thin,color=black!25] (N3_12) -- (M3_12);
\draw [thin,color=black!25] (N3_13) -- (M3_13);
\draw [thin,color=black!25] (N3_14) -- (M3_14);
\draw [thin,color=black!25] (N3_15) -- (M3_15);

\node (P3_10) at (-6,-4) [place] {};
\node (P3_11) at (-4.5,-4) [place1] {};
\node (P3_12) at (-3,-4) [place1] {};

\draw [thin,color=black!25] (P3_10) -- (P3_11);
\draw [thin,color=black!25] (P3_11) -- (P3_12);

\node (P3_13) at (-4.5,-5.5) [place1] {};
\node (P3_14) at (-3,-5.5) [place1] {};

\draw [thin,color=black!25] (P3_11) -- (P3_13);
\draw [thin,color=black!25] (P3_12) -- (P3_14);

\draw [thin,color=black!25] (P3_13) -- (P3_14);

\node (P3_15) at (-3,-7) [place1] {};

\draw [thin,color=black!25] (P3_14) -- (P3_15);

\draw [thin,color=black!25] (M3_10) -- (P3_10);
\draw [thin,color=black!25] (M3_11) -- (P3_11);
\draw [thin,color=black!25] (M3_12) -- (P3_12);
\draw [thin,color=black!25] (M3_13) -- (P3_13);
\draw [thin,color=black!25] (M3_14) -- (P3_14);
\draw [thin,color=black!25] (M3_15) -- (P3_15);

\node (Q3_13) at (-5.5,-6) [place1] {};
\node (Q3_14) at (-4,-6) [place1] {};

\draw [thin,color=black!25] (Q3_13) -- (Q3_14);

\node (Q3_15) at (-4,-7.5) [place1] {};

\draw [thin,color=black!25] (Q3_14) -- (Q3_15);

\draw [thin,color=black!25] (P3_13) -- (Q3_13);
\draw [thin,color=black!25] (P3_14) -- (Q3_14);
\draw [thin,color=black!25] (P3_15) -- (Q3_15);

\node (R3_15) at (-5,-8) [place1] {};

\draw [thin,color=black!25] (Q3_15) -- (R3_15);

\draw [dashed] (-1,-1.5) circle (0.4);
\draw [dashed] (-2,-2) circle (0.4);

\draw [dotted] (-0.6,-1.5) -- (-1.6,-2);
\draw [dotted] (-1,-1.1) -- (-2,-1.6);
\draw [dotted] (-1,-1.9) -- (-2,-2.4);

\draw [dashed] (-2.5,-4.5) circle (0.4);
\draw [dashed] (-3.5,-5) circle (0.4);

\draw [dotted] (-2.5,-4.1) -- (-3.5,-4.6);
\draw [dotted] (-2.5,-4.9) -- (-3.5,-5.4);
\draw [dotted] (-2.1,-4.5) -- (-3.1,-5);

\draw [->,thick,decorate,decoration={snake,amplitude=.3mm,segment length=2mm}] (M3_13) -- node[left]{\textbf{!}} (M3_11); 

\end{tikzpicture}
}
\quad
\subfloat[][$\widehat{B}=\{(x^2,xy^2,xz^2,y^4)_4\}$]{\label{fig:ExampleSubstMonom_3}
\begin{tikzpicture}[scale=0.65]
\tikzstyle{quotient}=[]
\tikzstyle{place}=[circle,draw=black,fill=black,inner sep=1pt]
\tikzstyle{place1}=[circle,draw=black,thick,inner sep=1pt]
\tikzstyle{place2}=[circle,draw=black!10,fill=black!30,inner sep=1pt]

\node (N2_1) at (1,0) [place] {};
\node (N2_2) at (2.5,0) [place] {};
\node (N2_3) at (4,0) [place] {};
\node (N2_4) at (5.5,0) [place] {};
\node (N2_5) at (7,0) [place] {};

\draw [thin,color=black!25] (N2_1) -- (N2_2);
\draw [thin,color=black!25] (N2_2) -- (N2_3);
\draw [thin,color=black!25] (N2_3) -- (N2_4);
\draw [thin,color=black!25] (N2_4) -- (N2_5);

\node (N2_6) at (2.5,-1.5) [place] {};
\node (N2_7) at (4,-1.5) [place] {};
\node (N2_8) at (5.5,-1.5) [place] {};
\node (N2_9) at (7,-1.5) [place1] {};

\draw [thin,color=black!25] (N2_2) -- (N2_6);
\draw [thin,color=black!25] (N2_3) -- (N2_7);
\draw [thin,color=black!25] (N2_4) -- (N2_8);
\draw [thin,color=black!25] (N2_5) -- (N2_9);

\draw [thin,color=black!25] (N2_6) -- (N2_7);
\draw [thin,color=black!25] (N2_7) -- (N2_8);
\draw [thin,color=black!25] (N2_8) -- (N2_9);

\node (N2_10) at (4,-3) [place] {};
\node (N2_11) at (5.5,-3) [place] {};
\node (N2_12) at (7,-3) [place1] {};

\draw [thin,color=black!25] (N2_7) -- (N2_10);
\draw [thin,color=black!25] (N2_8) -- (N2_11);
\draw [thin,color=black!25] (N2_9) -- (N2_12);

\draw [thin,color=black!25] (N2_10) -- (N2_11);
\draw [thin,color=black!25] (N2_11) -- (N2_12);

\node (N2_13) at (5.5,-4.5) [place] {};
\node (N2_14) at (7,-4.5) [place1] {};

\draw [thin,color=black!25] (N2_11) -- (N2_13);
\draw [thin,color=black!25] (N2_12) -- (N2_14);

\draw [thin,color=black!25] (N2_13) -- (N2_14);

\node (N2_15) at (7,-6) [place1] {};

\draw [thin,color=black!25] (N2_14) -- (N2_15);

\node (M2_6) at (1.5,-2) [place] {};
\node (M2_7) at (3,-2) [place] {};
\node (M2_8) at (4.5,-2) [place] {};
\node (M2_9) at (6,-2) [place1] {};

\draw [thin,color=black!25] (M2_6) -- (M2_7);
\draw [thin,color=black!25] (M2_7) -- (M2_8);
\draw [thin,color=black!25] (M2_8) -- (M2_9);

\node (M2_10) at (3,-3.5) [place] {};
\node (M2_11) at (4.5,-3.5) [place] {};
\node (M2_12) at (6,-3.5) [place1] {};

\draw [thin,color=black!25] (M2_7) -- (M2_10);
\draw [thin,color=black!25] (M2_8) -- (M2_11);
\draw [thin,color=black!25] (M2_9) -- (M2_12);

\draw [thin,color=black!25] (M2_10) -- (M2_11);
\draw [thin,color=black!25] (M2_11) -- (M2_12);

\node (M2_13) at (4.5,-5) [place] {};
\node (M2_14) at (6,-5) [place1] {};

\draw [thin,color=black!25] (M2_11) -- (M2_13);
\draw [thin,color=black!25] (M2_12) -- (M2_14);

\draw [thin,color=black!25] (M2_13) -- (M2_14);

\node (M2_15) at (6,-6.5) [place1] {};

\draw [thin,color=black!25] (M2_14) -- (M2_15);

\draw [thin,color=black!25] (N2_6) -- (M2_6);
\draw [thin,color=black!25] (N2_7) -- (M2_7);
\draw [thin,color=black!25] (N2_8) -- (M2_8);
\draw [thin,color=black!25] (N2_9) -- (M2_9);
\draw [thin,color=black!25] (N2_10) -- (M2_10);
\draw [thin,color=black!25] (N2_11) -- (M2_11);
\draw [thin,color=black!25] (N2_12) -- (M2_12);
\draw [thin,color=black!25] (N2_13) -- (M2_13);
\draw [thin,color=black!25] (N2_14) -- (M2_14);
\draw [thin,color=black!25] (N2_15) -- (M2_15);

\node (P2_10) at (2,-4) [place] {};
\node (P2_11) at (3.5,-4) [place] {};
\node (P2_12) at (5,-4) [place1] {};

\draw [thin,color=black!25] (P2_10) -- (P2_11);
\draw [thin,color=black!25] (P2_11) -- (P2_12);

\node (P2_13) at (3.5,-5.5) [place1] {};
\node (P2_14) at (5,-5.5) [place1] {};

\draw [thin,color=black!25] (P2_11) -- (P2_13);
\draw [thin,color=black!25] (P2_12) -- (P2_14);

\draw [thin,color=black!25] (P2_13) -- (P2_14);

\node (P2_15) at (5,-7) [place1] {};

\draw [thin,color=black!25] (P2_14) -- (P2_15);

\draw [thin,color=black!25] (M2_10) -- (P2_10);
\draw [thin,color=black!25] (M2_11) -- (P2_11);
\draw [thin,color=black!25] (M2_12) -- (P2_12);
\draw [thin,color=black!25] (M2_13) -- (P2_13);
\draw [thin,color=black!25] (M2_14) -- (P2_14);
\draw [thin,color=black!25] (M2_15) -- (P2_15);

\node (Q2_13) at (2.5,-6) [place1] {};
\node (Q2_14) at (4,-6) [place1] {};

\draw [thin,color=black!25] (Q2_13) -- (Q2_14);

\node (Q2_15) at (4,-7.5) [place1] {};

\draw [thin,color=black!25] (Q2_14) -- (Q2_15);

\draw [thin,color=black!25] (P2_13) -- (Q2_13);
\draw [thin,color=black!25] (P2_14) -- (Q2_14);
\draw [thin,color=black!25] (P2_15) -- (Q2_15);

\node (R2_15) at (3,-8) [place1] {};

\draw [thin,color=black!25] (Q2_15) -- (R2_15);

\draw [dashed] (7,-1.5) circle (0.4);
\draw [dashed] (6,-2) circle (0.4);

\draw [dotted] (7.4,-1.5) -- (6.4,-2);
\draw [dotted] (7,-1.1) -- (6,-1.6);
\draw [dotted] (7,-1.9) -- (6,-2.4);

\draw [dashed] (5.5,-4.5) circle (0.4);
\draw [dashed] (4.5,-5) circle (0.4);

\draw [dotted] (5.5,-4.1) -- (4.5,-4.6);
\draw [dotted] (5.5,-4.9) -- (4.5,-5.4);
\draw [dotted] (5.9,-4.5) -- (4.9,-5);
\end{tikzpicture}
}
\end{center}
\caption{\label{fig:ExampleSubstMonom} Some examples of substitutions in the poset $\mathcal{P}(3,4)$ breaking and preserving the Borel condition (see Example \ref{ex:breakPreserveBorel}). The monomials inside the Borel sets are the black points.}
\end{figure}
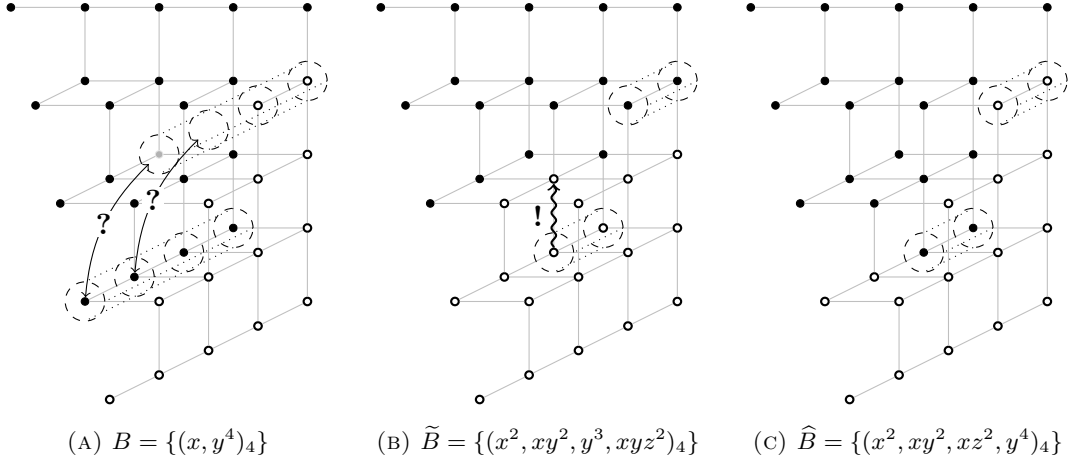 
\end{example}

Once determined how to swap monomials in the general case maintaining the Borel condition, we have to construct the deformation. We start with a technical lemma.

\begin{lemma}\label{lem:BorderMonomials}
 Let $B \subset \poset$ a Borel set, $x^\alpha$ a minimal monomial in $\mathfrak{m}(B_j)$ and $\mathcal{F}_\alpha$ the set of decreasing moves as in Definition \ref{def:decreasingMoves}. Moreover let us denote with $C$ the set $B\setminus\mathcal{F}_\alpha(x^\alpha)$. Then:
\begin{enumerate}
 \item\label{it:BorderMonomials_i} the monomials $x_i\, F(x^\alpha),\ i > j,\ F \in \mathcal{F}_\alpha$, belong to the ideal generated by $C$;
 \item\label{it:BorderMonomials_ii} the monomials $x_i\, F(x^\alpha),\ i \leqslant j,\ F \in \mathcal{F}_\alpha$, do not belong to the ideal generated by $C$. Furthermore a basis of $\langle x_j,\ldots,x_0 \rangle \cdot \langle \mathcal{F}_\alpha(x^\alpha)\rangle$ is given by
\begin{equation}
 \left\{x_i F(x^\alpha)\ \vert\ F \in \mathcal{F}_\alpha \text{ and } i = 0,\ldots,\min F(x^\alpha)\right\}.
\end{equation}
\end{enumerate}
\begin{proof}
\emph{(\ref{it:BorderMonomials_i})} Let $k = \min F(x^\alpha)$. $k \leqslant j < i$, then $x_i\, F(x^\alpha) = \up_{i-1} \cdot \ldots \cdot \up_{k} \big(F(x^\alpha)\big) \, x_k$ and $\up_{i-1} \cdot \ldots \cdot \up_{k} \big(F(x^\alpha)\big)$ belongs to $C$ by construction.

\emph{(\ref{it:BorderMonomials_ii})} Let us start considering $x^\alpha$. We could write for any $l \geqslant j = \min x^\alpha \geqslant i$ such that $x_l \mid x^\alpha$
\[
 x_i \, x^\alpha = \down_{i-1} \cdot \ldots \cdot \down_l (x^\alpha) x_l
\]
and $\down_{i-1} \cdot \ldots \cdot \down_l (x^\alpha)$ does not belong to $C$ for the minimality of $x^\alpha$, in fact if $\down_{i-1} \cdot \ldots \cdot \down_l \in \mathcal{F}_\alpha$ then it will belong to $\mathcal{F}_\alpha(x^\alpha)$, otherwise it will not belong to $B$.

Finally let us consider the generic $x^\gamma = F(x^\alpha)$ and let $k = \min x^\gamma$.
\begin{itemize}
 \item for $k < i \leqslant j$, $x_i \, x^\gamma = \up_{i-1}\cdot\ldots\cdot\up_{k} (x^\gamma) x_k$. But $\up_{i-1}\cdot\ldots\cdot\up_{k} (x^\gamma)$ has to belong to $\mathcal{F}_\alpha(x^\alpha)$ and $k \leqslant \min \big(\up_{i-1}\cdot\ldots\cdot\up_{k} (x^\gamma)\big)$.
 \item for $0 \leqslant i \leqslant k$, as done before for every $l \geqslant k$ such that $x_l \mid x^\gamma$, $x_i \, x^\gamma = \down_{i-1}\cdot\ldots\cdot\down_{l} (x^\gamma) x_l$ and $\down_{i-1}\cdot\ldots\cdot\down_{l} (x^\gamma) \notin B$. \qedhere
\end{itemize}
\end{proof}
\end{lemma}

\begin{theorem}\label{th:flatGeneralDim}
Let $\{I_r\}\subset\poset$ be a Borel set defined by the ideal $I$ with Hilbert polynomial $q(t)$. Moreover let us suppose that there exist two monomials $x^\alpha \in \mathfrak{m}\{I_r\}_j$ and $x^\beta \in \mathfrak{M}\{N(I)_r\}_j$, such that the set of decreasing moves $\mathcal{F}_\alpha$, as in Definition \ref{def:decreasingMoves}, is Borel-consistent w.r.t. $x^\beta$. The ideal 
\begin{equation}\label{eq:familyIdeals}
 \mathcal{I} = \left\langle \left\{I_r\right\}\setminus\mathcal{F}_\alpha(x^\alpha)\cup\{y_0 \, F(x^\alpha) + y_1 \, F(x^\beta) \ \vert\ F \in \mathcal{F}_\alpha\}\right\rangle
\end{equation}
has Hilbert polynomial $q(t)$, for all $[y_0 : y_1] \in \PP^1$.
\begin{proof}
Firstly by Proposition \ref{prop:BorelSetHP} the polynomial $p(t) = \binom{n+t}{t} - q(t)$ has the property $\left\vert \{N(I)_r\}_i \right\vert = \Delta^i p(r)$. Then Lemma \ref{lem:PreservingBorel} guarantees that also $\mathcal{I}\vert_{[0:1]} = \big\langle \{I_r\}\setminus\mathcal{F}_\alpha(x^\alpha)\cup\mathcal{F}_\alpha(x^\beta)\big\rangle$ has Hilbert polynomial $q(t)$.

As already done in the proof of Theorem \ref{th:flatZeroDim}, we use the Gotzmann Persistence Theorem (\ref{th:Gotzmann}). Starting from $\dim_k \mathcal{I}_r = q(r)$, we will prove $\dim_k \mathcal{I}_{r+1} = q(r+1)$.

Let us suppose $y_0 \neq 0$ and let us denote by $J$ the ideal $\mathcal{I}\vert_{[1:y_1]}$. We want to show that $ \dim_k J_{r+1} = \dim_k I_{r+1}$.
Let $C = B \setminus \mathcal{F}_\alpha(x^\alpha)$ and let $\mathcal{G}$ be the standard monomial basis of $I_{r+1}$. By Lemma \ref{lem:BorderMonomials}, we know that $\mathcal{G}$ can be split in two parts: $\mathcal{G}'$ of the monomials generated by the monomials in $C$ and $\mathcal{G}''$ of the monomials in $\langle x_j,\ldots,x_0 \rangle \cdot \langle \mathcal{F}_\alpha(x^\alpha)\rangle$.
Always by Lemma \ref{lem:BorderMonomials} we know that $\mathcal{G}'' = \left\{x_i F(x^\alpha)\ \vert\ F \in \mathcal{F}_\alpha \text{ and } i = 0,\ldots,\min F(x^\alpha)\right\}$ and we claim that 
\[
\widetilde{\mathcal{G}} = \mathcal{G}' \cup \left\{x_i F(x^\alpha) + y_1\, x_i F(x^\beta)\ \Big\vert\ F \in \mathcal{F}_\alpha \text{ and } i = 0,\ldots,\min F(x^\alpha)\right\}
\]
is a basis for $J_{r+1}$. To prove this claim we have only to check that the every polynomial $x_i \left(F(x^\alpha) + y_1\, F(x^\beta)\right)$ can be decomposed over $\widetilde{\mathcal{G}}$. Let $l = \min F(x^\alpha) = \min F(x^\beta)$.
\begin{itemize}
 \item For $i > j \geqslant l$
\[
 x_i\, F(x^\alpha) + y_1\, x_i\, F(x^\beta) = \up_{i-1}\cdot\ldots\cdot\up_{l} F(x^\alpha)\, x_{l} + y_1\, \up_{i-1}\cdot\ldots\cdot\up_{l} F(x^\beta)\, x_l
\]
and both $\up_{i-1}\cdot\dots\cdot\up_{l} F(x^\alpha)$ and $\up_{i-1}\cdot\dots\cdot\up_{l} F(x^\beta)$ belong to $C$ because of the hypothesis of minimality of $x^\alpha$ and maximality in $N_j$ of $x^\beta$.
\item For $0 \leqslant i \leqslant l$, $x_i\, F(x^\alpha) + y_1\, x_i \, F(x^\beta)$ belongs to $\widetilde{\mathcal{G}}$ by construction.
\item For $l < i \leqslant j$,
\[
 x_i\, F(x^\alpha) + y_1\, x_i\, F(x^\beta) = \up_{i-1}\cdot\ldots\cdot\up_{l} F(x^\alpha)\, x_{l} + y_1\, \up_{i-1}\cdot\ldots\cdot\up_{l} F(x^\beta)\, x_l
\]
and $\up_{i-1}\cdot\ldots\cdot\up_{l} F(x^\alpha)$ is in $\mathcal{F}_\alpha(x^\alpha)$, that is $\up_{i-1}\cdot\ldots\cdot\up_{l}\cdot F = G \in \mathcal{F}_\alpha$, so
\[
 x_i\, F(x^\alpha) + y_1\, x_i\, F(x^\beta) = x_l \, G(x^\alpha) + y_1\, x_l \, G(x^\beta),\quad G(x^\alpha) + y_1\, G(x^\beta) \in \mathcal{G}.
\]
\end{itemize}

In the case $y_1 \neq 0$, it suffices to start with the ideal $\mathcal{I}\vert_{[0:1]}$, that is exchanging the role of $x^\alpha$ and $x^\beta$, and to repeat the construction considering $\mathcal{F}_\beta\ (= \mathcal{F}_\alpha)$.
\end{proof}
\end{theorem}

\begin{theorem}\label{th:curveGeneralDim}
Let $\hilb^n_{p(t)}$ be the Hilbert scheme parametrizing subschemes of $\PP^n$ with Hilbert polynomial $p(t)$, whose Gotzmann number is $r$. Let $I$ and $J$ be two Borel-fixed ideals in $\hilbp$ and let $x^\alpha\in\mathfrak{m}\{I_r\}_i$, $x^\beta \in \mathfrak{M}\{N(I)_r\}_i$ and $\mathcal{F}_\alpha$ be as in Definition \ref{def:decreasingMoves}, such that $\mathcal{F}_\alpha$ is Borel-consistent w.r.t. $x^\beta$. If $\{J_r\} = \{I_r\} \cup \mathcal{F}_\alpha(x^\beta) \setminus\mathcal{F}_\alpha(x^\alpha)$, then there is a rational curve $\mathcal{C}: \PP^1 \rightarrow \hilbp$ going through the points of $\hilb^n_{p(t)}$ defined by $I$ and $J$.
\begin{proof}
Let us consider the morphism of rings
\[
 \widehat{\phi}: \K[y_0,y_1] \rightarrow \K[x][y_0,y_1]/\mathcal{I},\qquad \mathcal{I} = \left\langle \{I_r\}\setminus\mathcal{F}_\alpha(x^\alpha)\cup\{y_0\, F(x^\alpha) + y_1 \, F(x^\beta)\ \vert\ F \in \mathcal{F}_\alpha\}\right\rangle,
\]
 and the associated family $\phi$ of subschemes of $\PP^n$
\begin{center}
\begin{tikzpicture}
\node (a) at (0,1.5) {$\mathcal{Y}$};
\node (b) at (3,1.5) {$\PP^n \times \PP^1$};
\node (c) at (3,0) {$\PP^1$};
\draw [right hook->] (a) -- (b);
\draw [->] (a) -- node[below]{$\phi$} (c);
\draw [->] (b) -- (c);
\end{tikzpicture}
\end{center}
where $\mathcal{Y} = \Proj \K[x][y_0,y_1]/\mathcal{I}$. This family is flat by Theorem \ref{th:curveGeneralDim}, so there exists a map $\mathcal{C}: \PP^1 \rightarrow \hilb_{p(t)}^n$ that makes commutative the diagram:
\begin{center}
\begin{tikzpicture}
\node (a) at (0,1.5) {$\mathcal{Y}$};
\node (b) at (2.5,1.5) {$\mathcal{X}$};
\node (c) at (0,0) {$\PP^1$};
\node (d) at (2.5,0) {$\hilb_{p(t)}^n$};
\node at (1.25,0.75) {$\circlearrowleft$};
\draw [->] (a) -- (b);
\draw [->] (a) -- node[left]{$\phi$} (c);
\draw [->] (b) -- node[right]{$\h$} (d);
\draw [->] (c) -- node[above]{$\mathcal{C}$} (d);
\end{tikzpicture}
\end{center}
Finally, by construction the base extensions $\Proj \K[y_0,y_1]/(y_1) \rightarrow \PP^1$ and $\Proj \K[y_0,y_1]/(y_0) \rightarrow \PP^1$ specify the points $\h(\Proj \K[x]/I)$ and $\h(\Proj \K[x]/J)$.
\end{proof}
\end{theorem}

\begin{remark}
A rational deformation between two Borel-fixed ideals $I$ and $\widetilde{I}$ corresponds also to an edge connecting the vertices $I$ and $\widetilde{I}$ in the \emph{graph of monomial ideals} introduced in \cite{AltSturm}. The vertices of this graph are the monomial ideals in $\K[x]$ and two ideals $I,\widetilde{I}$ are connected by an edge if there exists an ideal $J$ such that the set of all its initial ideals (w.r.t. all term ordering) is $\{I,\widetilde{I}\}$.

Let $\mathcal{I}$ be the family of ideals described in \eqref{eq:familyIdeals} defining the deformation from $I$ to $\widetilde{I}$. By construction both $x^\alpha \geq_B x^\beta$ and $x^\alpha \leq x^\beta$ do not hold, so varying the term ordering $\preceq$ both relation $x^\alpha \succ x^\beta$ and $x^\alpha \prec x^\beta$ can be obtained. 

Let $k = \min x^\alpha = \min x^\beta$ and let $s = \max \{\lambda_k \ \vert\ (\down_1)^{\lambda_1}\cdot\ldots\cdot(\down_k)^{\lambda_k} = F \in \mathcal{F}_{\alpha}\}$. We know that $x_k^s \mid \text{gcd}(x^\alpha,x^\beta)$, then $x^\alpha = x^{\alpha'}\, x_k^s$, $x^\beta = x^{\beta'}\, x_k^s$ and $F(x^\alpha) = x^{\alpha'}\, F(x_k^s)$, $F(x^\beta) = x^{\beta'}\, F(x_k^s)$. Finally
\[
 x^\alpha = x^{\alpha'}\, x^s_k \succ x^{\beta'}\, x_k^s = x^\beta \quad \Rightarrow\quad x^{\alpha'} \succ x^{\beta'} \quad\Rightarrow\quad F(x^\alpha) = x^{\alpha'}\, F(x_k^s)\succ x^{\beta'}\, F(x_k^s) = F(x^\beta),
\]
$\forall\ F \in \mathcal{F}_\alpha$, and vice versa if $x^\alpha \prec x^\beta$.
\end{remark}

\begin{algorithm}
\caption{\label{alg:allDeformations} Pseudo-code description of the algorithm computing all the possible rational deformations of a given Borel-fixed ideal $I$. A trial version of this algorithm is available at the web page \href{http://www.dm.unito.it/dottorato/dottorandi/lella/HSC/deformations.html}{www.dm.unito.it/dottorato/dottorandi/lella/HSC/deformations.html}}.
\begin{algorithmic}
\Require $I\subset\K[x]$, Borel-fixed ideal.
\Ensure set of Borel-fixed ideals obtainable from $I$ by a rational deformation.
\Procedure{RationalDeformations}{$I$}
\State $r \leftarrow$ Gotzmann number of the Hilbert polynomial of $\K[x]/I$;
\State $S \leftarrow \emptyset$;
\For{$i=0,\ldots,n-1$}
\ForAll{$x^\alpha \in \mathfrak{m}\{I_r\}_i,\, x^\beta \in \mathfrak{M}\{N(I)_r\}_i$ s.t $\nexists\ k > i,\ \down_k(x^\alpha)=x^\beta$}
\State$\mathcal{F}_{\alpha} \leftarrow$ set of admissible decreasing moves w.r.t. $x^\alpha$;
\If{$\mathcal{F}_{\alpha}$ is Borel-consistent w.r.t. $x^\beta$}
 \State $S \leftarrow S \cup \big\{\big\langle\{I_r\}\cup\mathcal{F}_{\alpha}(x^\beta)\setminus\mathcal{F}_{\alpha}(x^\alpha)\big\rangle\big\}$;
\EndIf
\EndFor
\EndFor
\State \Return $S$;
\EndProcedure
\end{algorithmic}
\end{algorithm}

\begin{example}\label{ex:multipleDef}
 Let us simulate an execution of Algorithm \ref{alg:allDeformations} on the ideal $I = (x_3^2,x_3 x_2^2, x_3 x_2 x_1,$ $x_2^4,x_2^3 x_1, x_2^2 x_1^2)_{\geqslant 8} \subset \K[x_0,x_1,x_2,x_3]$. The corresponding subscheme $\K[x_0,x_1,x_2,x_3]/I$ has Hilbert polynomial $p(t)=3t+5$, whose Gotzmann number is $8$.
\begin{enumerate}
 \item[$i=0$.] $\mathfrak{m}\{I_8\}_0 = \{x_3^2 x_0^6,x_3 x_2 x_1 x_0^5,x_2^2 x_1^2 x_0^4\}$ and $\mathfrak{M}\{N(I)_8\}_0 = \{x_3 x_2 x_0^6,x_2^3 x_0^5\}$.
  \begin{itemize}
  \item $x_3^2 x_0^6$ and $x_3 x_2 x_0^6$. $\down_3(x_3^2 x_0^6) = x_3 x_2 x_0^6$.
  \item $x_3^2 x_0^6$ and $x_2^3 x_0^5$. $\mathcal{F} = \{\mathrm{id}\}$ is Borel-consistent w.r.t. $x_2^3 x_0^5$:
\[
 J_1 = \left\langle\{I_8\}\setminus\{x_3^2 x_0^6\}\cup\{x_2^3 x_0^5\}\right\rangle =(x_3^3,x_3^2 x_2, x_3 x_2^2, x_2^3,x_3^2 x_1, x_3 x_2 x_1, x_2^2 x_1^2)_{\geqslant 8}.
\]
  \item $x_3 x_2 x_1 x_0^5$ and $x_3 x_2 x_0^6$. $\down_1(x_3 x_2 x_1 x_0^5) = x_3 x_2 x_0^6$.
  \item $x_3 x_2 x_1 x_0^5$ and $x_2^3 x_0^5$. $\mathcal{F} = \{\mathrm{id}\}$ is Borel-consistent w.r.t. $x_2^3 x_0^5$:
\[
 J_2 = \left\langle\{I_8\}\setminus\{x_3 x_2 x_1 x_0^5\}\cup\{x_2^3 x_0^5\}\right\rangle = (x_3^2,x_3 x_2^2,x_2^3,x_3 x_2 x_1^2,x_2^2 x_1^2)_{\geqslant 8}.
\]
  \item $x_2^2 x_1^2 x_0^4$ and $x_3 x_2 x_0^6$. $\mathcal{F} = \{\mathrm{id}\}$ is Borel-consistent w.r.t. $x_3 x_2 x_0^6$:
\[
 J_3 = \left\langle\{I_8\}\setminus\{x_2^2 x_1^2 x_0^4\}\cup\{x_3 x_2 x_0^6\}\right\rangle = (x_3^2,x_3 x_2,x_2^4,x_2^3 x_1, x_2^2 x_1^3)_{\geqslant 8}. 
\]
  \item $x_2^2 x_1^2 x_0^4$ and $x_2^3 x_0^5$. $\mathcal{F} = \{\mathrm{id}\}$ is Borel-consistent w.r.t. $x_2^3 x_0^5$:
\[
 J_4 = \left\langle\{I_8\}\setminus\{x_2^2 x_1^2 x_0^4\}\cup\{x_2^3 x_0^5\}\right\rangle = (x_3^2,x_3 x_2^2,x_2^3,x_3 x_2 x_1,x_2^2,x_1^3)_{\geqslant 8}. 
\]
  \end{itemize}
 \item[$i=1$.] $\mathfrak{m}\{I_8\}_1 = \{x_2^2 x_1^6\}$ and $\mathfrak{M}\{N(I)_8\}_1 = \{x_3 x_1^7\}$.
  \begin{itemize}
  \item $x_2^2 x_1^6$ and $x_3 x_1^7$. $\mathcal{F} = \{\mathrm{id},\down_1,(\down_1)^2,(\down_1)^3,(\down_1)^4\}$ is Borel-consistent w.r.t. $x_3 x_1^7$:
\[
 J_5 = \left\langle\{I_8\}\setminus\mathcal{F}(x_2^2 x_1^6)\cup\mathcal{F}(x_3 x_1^7)\right\rangle = (x_3^2,x_3 x_2^2,x_3 x_2 x_1,x_2^4,x_2^3 x_1,x_3 x_1^3)_{\geqslant 8}.
\]
  \end{itemize}
 \item[$i=2.$] $\mathfrak{m}\{I_8\}_2 = \{x_2^8\}$ and $\mathfrak{M}\{N(I)_8\}_2 = \emptyset$.
\end{enumerate}
\end{example}

\section{The connectedness of the Hilbert scheme}\label{sec:connected}
In this section, we are going to study consecutive deformations of Borel ideals and we want to have a control on the \lq\lq direction\rq\rq\ toward we move, in order to obtain a technique similar to the one introduced by Peeva and Stillman \cite{PS}. Those deformations are affine and based on Gr\"obner basis tool, but anyhow the idea is to exchange monomials belonging to the ideal with monomials not belonging. The goal is to determine a sequence of deformations leading from any Borel-fixed ideal to the lexicographic ideal, so the choice of the monomials to exchange is governed by the $\mathtt{DegLex}$ term ordering.

Therefore we start slightly modifying Algorithm \ref{alg:allDeformations}, adding a term ordering $\preceq$, refinement of the Borel partial order $\leq_B$, that we will use to choose in a unique way the monomials to swap.

\begin{algorithm}[H]
\caption{\label{alg:deformation} How to find a rational deformation w.r.t. a fixed term ordering. A trial version of this algorithm is available at the web page\\ \href{http://www.dm.unito.it/dottorato/dottorandi/lella/HSC/TOdeformation.html}{www.dm.unito.it/dottorato/dottorandi/lella/HSC/TOdeformation.html}}
\begin{algorithmic}
\Require  $I\subset\K[x]$, Borel-fixed ideal.
\Require $\preceq$, term ordering, refinement of the Borel partial order $\leq_B$.
\Ensure the ideal $J$, deformation of $I$ w.r.t. $\preceq$.
\Procedure{RationalDeformation}{$I$,$\preceq$}
\State $r \leftarrow$ Gotzmann number of the Hilbert polynomial of $\K[x]/I$;
\For{$i=0,\ldots,n-1$}
\If{$x^\alpha = \min_\preceq \mathfrak{m}\{I_r\}_i \prec x^\beta = \max_\preceq \mathfrak{M}\{N(I)_r\}_i$}
\State$\mathcal{F}_{\alpha} \leftarrow$ set of admissible decreasing moves w.r.t. $x^\alpha$;
\If{$\mathcal{F}_{\alpha}$ is Borel-consistent w.r.t. $x^\beta$}
 \State \Return $J = \langle\{I_r\}\cup\mathcal{F}_{\alpha}(x^\beta)\setminus\mathcal{F}_{\alpha}(x^\alpha)\big\rangle\big\}$;
\EndIf
\EndIf
\EndFor
\State \Return $I$;
\EndProcedure
\end{algorithmic}
\end{algorithm}

We remark that the condition $x^\alpha = \min_\preceq \mathfrak{m}\{I_r\}_i \prec x^\beta = \max_\preceq \mathfrak{M}\{N(I)_r\}_i$ guarantees that $x^\beta$ can not be obtained from $x^\alpha$ by an elementary move $\down_k$. In fact if such a move exists, $x^\alpha \geq_B x^\beta$ implies $x^\alpha \succ x^\beta$ for any term ordering $\preceq$. Moreover the deformation determined is unique because the algorithm returns the first one found. Algorithm \ref{alg:deformation} could return also the same ideal $I$ given as input, and this is a reasonable possibility, because for instance if we apply the algorithm on the lexsegment ideal with the $\mathtt{DegLex}$ term ordering, since $x^\alpha >_{\mathtt{DegLex}} x^\beta, \forall\ x^\alpha \in I_r,\ x^\beta \in N(I)_r$ by definition of lexsegment, the condition $x^\alpha = \min_\preceq \mathfrak{m}\{I_r\}_i <_{\mathtt{DegLex}} x^\beta = \max_\preceq \mathfrak{M}\{N(I)_r\}_i$ will be always false.

\begin{definition}\label{def:preceqdef}
 Given a Borel ideal $I$ and a term ordering $\preceq$, we say that
\begin{itemize}
 \item the ideal $J \neq I$ returned by $\textsc{RationalDeformation}(I,\preceq)$ (Algorithm \ref{alg:deformation}) is a $\preceq$-\emph{rational deformation} or simply a $\preceq$-\emph{deformation} of $I$;
 \item $I$ is a $\preceq$-\emph{endpoint} if $\textsc{RationalDeformation}(I,\preceq)$ returns the same ideal $I$.
\end{itemize}
\end{definition}

After having determined a method similar to the one proposed in \cite{PS}, we want to compare this two approaches, so we recall briefly the strategy and some notation of that paper. Given an ideal $I = I_{\geqslant r} \subset \K[x]$, Peeva and Stillman compute the monomials $x^\beta = \max_{\mathtt{DegLex}} N(I)_r$ (they call it \emph{first gap}) and $x^\alpha = \max_{\mathtt{DegLex}} \{x^\delta \in I_r\ \vert\ x^\delta <_{\mathtt{DegLex}} x^\beta\}$. If $I$ is the lexsegment ideal, the set $\{x^\delta \in I_r\ \vert\ x^\delta <_{\mathtt{DegLex}} x^\beta\}$ will be obviously empty. At this point they determine the set of monomials 
\[
T = \{x^{\alpha'} x^{\gamma}\in I_r \ \vert\ x^\alpha = x^{\alpha'}\cdot x_{\min x^\alpha}^{s},\ x^{\beta'} x^\gamma \notin I_r \text{ and } \max x^\gamma < \min x^{\beta'}\}   
\]
and they fix the monomial $x^{\alpha'}$ of minimal degree and the corresponding $x^{\beta'}$ of minimal degree. Let us denote them by $x^{\overline{\alpha}}$ and $x^{\overline{\beta}}$. They finally form the ideal
\[
 \widetilde{I} = \left\langle\{I_r\}\setminus\{x^{\overline{\alpha}} x^{\gamma} \in I_r\}\cup\{x^{\overline{\alpha}} x^{\gamma} - x^{\overline{\beta}} x^{\gamma} \ \vert\ x^{\overline{\alpha}} x^{\gamma} \in I_r\}\right\rangle,\quad \max x^\gamma \leqslant \min x^{\overline{\beta}}.
\]
The Borel ideal \lq\lq lex-closer\rq\rq\ to the lexsegment ideal is $\text{gin}_{\mathtt{Lex}}\big(\text{in}_{\mathtt{Lex}}(\widetilde{I})\big)$.

The computation of a generic initial ideal reveals an important difference between the two techniques: this choice of monomials involved in the substitution generally does not preserve the Borel condition, so that to restore it a computation of a gin could be needed, as the following example shows.

\begin{example}\label{ex:PSbreakingBorel}
 Let us consider the ideal $I = (x^4,x^3y,x^2 y^2,x y^3)_{\geqslant 7} \subset \K[x,y,z]$. The Hilbert polynomial is $p(t) = t+7$ with Gotzmann number 7. The first gap is $x^3z^4$ and the lex-greatest monomial in $I_7$, lex-smaller than the gap, is $x^2 y^5$.
\[
 T = \left\{x^2 y^2 z^3, x^2 y^3 z^2, x^2 y^4 z, x^2 y^5\right\}
\]
so that $x^{\overline{\alpha}} = x^2 y^2$ and $x^{\overline{\beta}} = x^3 z$. We construct the ideal 
\[
\widetilde{I} = \left\langle \{I_7\} \setminus\{x^2 y^2 z^3\}\cup\{x^2 y^2 z^3- x^3 z^4\} \right\rangle,
\]
and $\text{in}_{\mathtt{Lex}}(\widetilde{I}) = (x^3,xy^3)_{\geqslant 7}$ is not Borel, because $x y^3 z^3 \in \text{in}_{\mathtt{Lex}}(\widetilde{I})$  and $\up_y(x y^3 z^3) = x^2 y^2 z^3 \notin \text{in}_{\mathtt{Lex}}(\widetilde{I})$. Finally $\text{gin}_{\mathtt{Lex}}\big(\text{in}_{\mathtt{Lex}}(\widetilde{I})\big) = (x^3,x^2 y^2,x y^4)_{\geqslant 7}$ is again Borel-fixed.

Using Algorithm \ref{alg:deformation} to compute the $\mathtt{DegLex}$-deformation of $I$, we have
\[
\min_{\mathtt{DegLex}} \mathfrak{m}\{I_7\}_0 =  xy^3z^3 \prec_{\mathtt{DegLex}} x^3 z^4 = \max_{\mathtt{DegLex}} \mathfrak{M}\{N(I)_7\}_0,
\]
and since the only admissible decreasing move is the identity, the $\mathtt{DegLex}$-rational deformation of $I$ is directly the ideal $J = \left\langle \{I_7\}\cup\{x^3 z^4\}\setminus\{xy^3z^3\} \right\rangle = (x^3,x^2 y^2,x y^4)_{\geqslant 7}$.
\end{example}

The choice of the first gap is very similar to the choice of the maximal monomial in $N$, because in both cases we look for a \lq\lq greatest\rq\rq\ element. Instead the two techniques differ in the choice of the monomial inside the ideal: we look for minimal elements whereas Peeva and Stillman consider maximal elements lower than the gap. Hence we expect that in the cases in which the monomials chosen are the same, the deformation is almost equal, in the sense that the the monomial ideal we reach is the same.

\begin{example}\label{ex:PSsameRational}
We consider the ideal $I = (x_3^2,x_3 x_2,x_2^3,x_2^2 x_1)_{\geqslant 5} \subset \K[x_0,x_1,x_2,x_3]$. The Hilbert polynomial is $p(t) = 3t+2$ with Gotzmann number 5. The first gap is $x_3 x_1^4$ and the greatest monomial of the ideal smaller than it is $x_2^5$. Then
\[
 T = \left\{ x_2^2 x_1^3, x_2^2 x_1^2 x_0, x_2^2 x_1 x_0^2, x_2^3 x_1^2, x_2^3 x_1 x_0, x_2^3 x_0^2, x_2^4 x_1, x_2^4 x_0, x_2^5\right\},
\]
$x^{\overline{\alpha}} = x_2^2$ and $x^{\overline{\beta}} = x_3 x_1$, so that 
\[
 \widetilde{I} = \left\langle I_5 \setminus\left\{x_2^2 x_1^3, x_2^2 x_1^2 x_0, x_2^2 x_1 x_0^2\right\}\cup\left\{x_2^2 x_1^3 - x_3 x_1^4, x_2^2 x_1^2 x_0 - x_3 x_1^3 x_0, x_2^2 x_1 x_0^2 - x_3 x_1^2 x_0^2\right\} \right\rangle
\]
and $\text{in}_{\mathtt{Lex}}(\widetilde{I}) = \text{gin}_{\mathtt{Lex}}\big(\text{in}_{\mathtt{Lex}}(\widetilde{I})\big) = (x_3^2, x_3 x_2, x_3 x_1^2,x_2^3)_{\geqslant 5}$.

Applying Algorithm \ref{alg:deformation} on $I$ with the $\mathtt{DegLex}$ term ordering, we have
\[
 \begin{split}
\min_{\mathtt{DegLex}} \mathfrak{m}\{I_5\}_0 = x_2^2 x_1 x_0^2 &\succ_{\mathtt{DegLex}} x_2^2 x_0^3 = \max_{\mathtt{DegLex}} \mathfrak{M}\{N(I)_5\}_0,\\
 \min_{\mathtt{DegLex}} \mathfrak{m}\{I_5\}_1 = x_2^2 x_1^3 &\prec_{\mathtt{DegLex}} x_3 x_1^4 = \max_{\mathtt{DegLex}} \mathfrak{M}\{N(I)_5\}_1. 
 \end{split}
\]
Then the algorithm determines $\mathcal{F} = \{\textrm{id},\down_1,(\down_1)^2\}$ that is Borel-consistent w.r.t. $x_3 x_1^4$, and swapping $\mathcal{F}(x_2^2 x_1^3) = \{x_2^2 x_1^3,x_2^2 x_1^2 x_0,x_2^2 x_1 x_0^2\}$ with $\mathcal{F}(x_3 x_1^4) = \{x_3 x_1^4,x_3 x_1^3 x_0,x_3 x_1^2 x_0^2\}$ we obtain again the ideal
\[
 J = (x_3^2, x_3 x_2, x_3 x_1^2,x_2^3)_{\geqslant 5} = \text{in}_{\mathtt{Lex}}(\widetilde{I}).
\]
\end{example}

Our next goal is to prove that also with our method it is possible to reach from any Borel-fixed ideal some special point on the Hilbert scheme, by a sequence of consecutive rational deformations. Of course the lexicographic ideal defines a special point, but we want to generalize the property to a wider class of ideals. Firstly we recall a definition given in \cite{CLMR}.

\begin{definition}\label{def:hilbSegment}\cite[Definition 3.7]{CLMR} Let $I \subset \K[x]$ be a Borel-fixed ideal and let $p(t)$ be the Hilbert polynomial of $\K[x]/I$ with Gotzmann number equal to $r$. $I$ is said to be a \emph{hilb-segment} if there exists a term ordering $\preceq$ such that $I_r$ contains the greatest monomials of degree $r$ w.r.t. $\preceq$, that is
\begin{equation}
 x^\alpha \succ x^\beta \qquad\forall\ x^\alpha \in I_r,\ \forall\ x^\beta \in N(I)_r.
\end{equation}
\end{definition}

Concretely, it is quite easy to determine whenever a Borel ideal $I$ could be a hilb-segment ideal or not. In fact, computed the Gotzmann number $r$, we can consider the Borel set $\{I_r\}$ with the minimal monomials $x^{\alpha_1},\ldots,x^{\alpha_u} \in \{I_r\}$ and the maximal ones $x^{\beta_1},\ldots,x^{\beta_s} \notin \{I_r\}$. We just look for an array $\omega = (\omega_n,\ldots,\omega_0) \in \mathbb{Z}^{n+1}$, $\omega_n > \ldots > \omega_0 > 0$, such that
\begin{equation*}
 \omega \cdot \alpha_i > \omega \cdot \beta_j, \qquad\forall\ i = 1,\ldots,u,\ j = 1,\ldots,s.
\end{equation*}
This is a system of linear inequalities that can be solved by applying the simplex algorithm, so that if a solution exists the term ordering $\preceq$ defined by the matrix
\begin{equation}\label{eq:matrixTO}
 \left[\begin{array}{c c c c c}
        1 & \ldots & 1 & \ldots & 1 \\
        \omega_n & \ldots & \omega_i & \ldots & \omega_0 \\
        0 & 1 & 0 & \ldots & 0 \\
        & & \ddots & & \\
        0 &\ldots & 0 & 1 & 0
       \end{array}
\right]
\end{equation}
refines the Borel partial order and makes $I$ a hilb-segment (a demo version of this algorithm is available at the web page \href{http://www.dm.unito.it/dottorato/dottorandi/lella/HSC/segment.html}{www.dm.unito.it/dottorato/dottorandi/lella/HSC/segment.html}).

\medskip

Now chosen a Hilbert scheme $\hilbp$ and a term ordering $\preceq$, we want to have an overview on all $\preceq$-rational deformations among Borel-fixed ideals belonging to $\hilbp$.

\begin{definition}
 Let $\hilbp$ be the Hilbert scheme parametrizing the subschemes of the projective space $\PP^n$ with Hilbert polynomial $p(t)$ and let $\preceq$ be any term ordering. We define the \emph{$\preceq$-deformation graph} of $\hilbp$ as the graph $(V,E)$, such that
\begin{itemize}
 \item the set $V$ of \emph{vertices} contains all Borel-fixed ideals in $\hilbp$;
 \item the set $E$ of \emph{edges} contains the $\preceq$-rational deformations connecting two Borel ideals.
\end{itemize}
The algorithm computing the $\preceq$-deformation graph is described in Algorithm \ref{alg:deformationGraph}.
\end{definition}
 
\begin{algorithm}[H]
\caption{\label{alg:deformationGraph} How to compute the $\preceq$-deformation graph. A trial version of this algorithm is available at the web page\\ \href{http://www.dm.unito.it/dottorato/dottorandi/lella/HSC/deformationGraph.html}{www.dm.unito.it/dottorato/dottorandi/lella/HSC/deformationGraph.html}.}
\begin{algorithmic}
\Require $\hilbp$, Hilbert scheme parametrizing subschemes of $\PP^n$ with Hilbert polynomial $p(t)$.
\Require $\preceq$, term ordering, refinement of the Borel partial order $\leq_B$.
\Ensure the $\preceq$-deformation graph $(V,E)$.
\Procedure{DeformationGraph}{$\hilbp$,$\preceq$}
\State $V\, \leftarrow$ Borel-fixed ideals in $\K[x]$ with Hilbert polynomial $p(t)$.
\State $E\, \leftarrow \emptyset$;
\ForAll{$I \in V$}
\State $J \leftarrow$ \Call{RationalDeformation}{$I$,$\preceq$};
\If{$I = J$}
\State $I \leftarrow$ $\preceq$-endpoint;
\Else
\State $E \leftarrow E \cup \{(I,J)\}$;
\EndIf
\EndFor
\State\Return $(V,E)$;
\EndProcedure
\end{algorithmic}
\end{algorithm}

\begin{proposition}\label{prop:connectedTree}
Let $\preceq$ be a term ordering such that the Hilbert scheme $\hilbp$ contains a Borel-fixed ideal $I$, hilb-segment w.r.t. $\preceq$. Then the $\preceq$-deformation graph of $\hilbp$ is a rooted tree, with $I$ as root.
\end{proposition}
We recall briefly what we mean by rooted tree. A \emph{tree} in connected graph $(V,E)$, such that $\vert E \vert = \vert G \vert -1$. A \emph{rooted tree} is a tree in which a fixed vertex (the root) determines a natural orientation of the edges, \lq\lq toward to\rq\rq\ and \lq\lq away from\rq\rq\ the root.
\begin{proof}
By definition $I$ is a $\preceq$-endpoint, so $I$ is the natural root because $I$ have not a deformation w.r.t. $\preceq$. To prove that the $\preceq$-deformation graph $(V,E)$ it is sufficient to show that any other Borel ideal $J \neq I$ has a $\preceq$-rational deformation.

Let $r$ be the Gotzmann number of $p(t)$. Let us suppose that for $j=s+1,\ldots,n-1$,
\[
x^\gamma \succ x^\delta,\quad \forall\ x^\gamma \in \{J_r\}_j,\ \forall\ x^\delta \in \{N(J)_r\}_j,
\]
and that there exists a monomial in $\{J_r\}_s$ smaller than a monomial in $\{N(J)_r\}_s$. Hence we know that
\[
 x^\beta = \max_{\preceq}\mathfrak{M}\{N(J)_r\}_s \succ \min_{\preceq}\mathfrak{m}\{J_r\}_s = x^\alpha.
\]
Let $\mathcal{F}_\alpha$ be the set of decreasing moves w.r.t. $x^\alpha$. If $\mathcal{F}_\alpha$ is Borel-consistent w.r.t. $x^\beta$, we finish because we are sure that Algorithm \ref{alg:deformation} applied on $J$ and $\preceq$ does not return the same ideal $J$. Then let us consider the case that $\mathcal{F}_\alpha$ is not Borel-consistent: we want to show that there must exist a deformation determined by a couple of monomials $x^{\alpha'}$ and $x^{\beta'}$ such that $\min x^{\alpha'} = \min x^{\beta'} < s$.
\begin{itemize}
 \item Let us start supposing that there exists some decreasing move in $\mathcal{F}_\alpha$ not admissible w.r.t. $x^\beta$. Let $\mathcal{G} \subset \mathcal{F}_\alpha$ be the set of the decreasing moves admissible on both monomials. Since $x^\alpha \in J$ does not belong to the hilb-segment $I$, also every monomial in $\mathcal{F}_\alpha$ does not belong to $I$. The monomials in $\mathcal{G}(x^\alpha)$ are replaced by the monomials in $\mathcal{G}(x^\beta)$ and the others in $\mathcal{F}_\alpha(x^\alpha)\setminus\mathcal{G}(x^\alpha)$ have to be replaced by monomials not obtained by decreasing moves from monomials in $\{J_r\}_s$, that is $\max_{\preceq}\mathfrak{M}\{N(J)_r\}_i \succ \min_{\preceq}\mathfrak{m}\{J_r\}_i,\ i < s$.
 \item Let us now consider the case of a set of decreasing moves $\mathcal{F}_\alpha$ admissible w.r.t. $x^\beta$ but such that for a certain $F \in \mathcal{G},\ \exists\ \up_k,\ k \geqslant s$, $\up_k \big(F(x^\beta)\big) \notin J$. Let $i = \min F(x^\alpha) = \min F(x^\beta) < s$. $\up_l F(x^\beta)$ can not be obtained by a monomial in $\{N(J)_r\}_s$ applying a composition of decreasing elementary moves $G$, because $G(x^\delta) = \up_l F(x^\beta)$ implies $x^\delta \geq_B x^\beta$ in contradiction with the hypothesis $x^\beta \in \mathfrak{M}\{N(J)_r\}_s$. Since $x^\alpha \prec x^\beta \Rightarrow F(x^\alpha) \prec F(x^\beta)$, $\max_\preceq \mathfrak{M}\{N(J)_r\}_i \succeq \up_k\big(F(x^\beta)\big) \succ F(x^\beta) \succ F(x^\alpha) \succeq \min_\preceq \mathfrak{m}\{J_r\}_i$.
\end{itemize}
In both cases, there exists $i < s$ and a couple of monomials $x^{\alpha'},x^{\beta'}$
\[
  x^{\beta'} = \max_{\preceq}\mathfrak{M}\{N(J)_r\}_i \succ \min_{\preceq}\mathfrak{m}\{J_r\}_i = x^{\alpha'}.
\]
We compute again $\mathcal{F}_{\alpha'}$ and we check if it is Borel-consistent w.r.t. $x^{\beta'}$: if not we repeat the reasoning and we look for monomials involving more variables. Finally, we are sure to find a Borel-consistent set of decreasing moves because if we reach the smallest variable $x_0$, the set contains only the identity move. 

The total order on the monomials guarantees that it is not possible to have cycles, because each deformation approaches to the hilb-segment.
\end{proof}

\begin{corollary}
 The Hilbert scheme $\hilbp$ is connected.
\begin{proof}
Let $I$ any ideal corresponding to a point on $\hilbp$. As usually through an affine deformation the point defined by $I$ can be connected to the point defined by the Borel-fixed ideal $\text{gin}(I)$.

Of course on $\hilbp$, there is the lexicographic point corresponding to the lexicographic ideal, that is hilb-segment w.r.t. $\mathtt{DegLex}$. By Proposition \ref{prop:connectedTree}, we know that the $\mathtt{DegLex}$-deformation graph is a connected rooted tree, so the point defined by any ideal $I$ can be connected to the lexicographic point by an initial affine deformation and a sequence of $\mathtt{DegLex}$-rational deformations.
\end{proof}
\end{corollary}

We underline that if $\preceq$ has no hilb-segment on $\hilbp$, the deformation graph could be not connected, as the following example shows.
 
\begin{example}
 Let us consider the Hilbert scheme $\hilb_{6t-5}^3$. There are 11 Borel-fixed ideals in $\K[x_0,x_1,x_2,x_3]$ with Hilbert polynomial $6t-5$:
\[
 \begin{array}{l c l}
  I_{1} = (x_3,x_2^7,x_2^6x_1^4)_{\geqslant 10}, && I_{7} = (x_3^2,x_3x_2^2,x_3x_2x_1,x_3x_1^2,x_2^7,x_2^6x_1)_{\geqslant 10}, \\
  I_{2} = (x_3,x_2^8,x_2^7x_1,x_2^6x_1^3)_{\geqslant 10}, && I_{8} = (x_3^2,x_3x_2,x_3x_1^4,x_2^6)_{\geqslant 10}, \\
  I_{3} = (x_3^2,x_3x_2,x_3x_1,x_2^7,x_2^6x_1^3)_{\geqslant 10}, && I_{9} = (x_3^2,x_3x_2^2,x_3x_2x_1,x_3x_1^3,x_2^6)_{\geqslant 10}, \\
  I_{4} = (x_3^2,x_3x_2,x_3x_1,x_2^8,x_2^7x_1,x_2^6x_1^2)_{\geqslant 10}, && I_{10} = (x_3^3,x_3^2x_2,x_3x_2^2,x_3^2x_1,x_3x_2x_1,x_3x_1^2,x_2^6)_{\geqslant 10}, \\
  I_{5} = (x_3^2,x_3x_2,x_2x_1^2,x_2^7,x_2^6x_1^2)_{\geqslant 10}, && I_{11} = (x_3^2,x_3x_2,x_2^5)_{\geqslant 10}, \\
  I_6 =(x_3^2,x_3x_2,x_3x_1^3,x_2^7,x_2^6x_1)_{\geqslant 10}. && 
 \end{array}
\]
Many of them are hilb-segments: for instance $I_5$ is a hilb-segment w.r.t. the term ordering $\preceq_5$ defined by the matrix \eqref{eq:matrixTO} with the array $\omega=(25,5,2,1)$ and the $\preceq_5$-deformation graph is a rooted tree (Figure \ref{fig:BorelGraph}\subref{fig:BorelGraph_1}). 

The\hfill ideal\hfill generates\hfill by\hfill the\hfill greatest\hfill $q(10) = 231$\hfill monomials\hfill in\hfill $\K[x_0,x_1,x_2,x_3]_{10}$\hfill w.r.t.\\ $\mathtt{RevLex}$ has constant Hilbert polynomial equal to $55$, so on $\hilb^3_{6t-5}$ there is not a hilb-segment w.r.t. $\mathtt{RevLex}$. Both $I_{10}$ and $I_{11}$ are $\mathtt{RevLex}$-endpoint, so that the $\mathtt{RevLex}$-deformation graph is not connected (Figure \ref{fig:BorelGraph}\subref{fig:BorelGraph_2}).
\begin{figure}[!ht] 
\begin{center}\captionsetup[subfloat]{width=5.5cm}
\subfloat[][The $\preceq_5$-deformation graph.]{\label{fig:BorelGraph_1}
\begin{tikzpicture}[>=latex,scale=0.5]
\tikzstyle{place}=[ellipse,draw=black!50,minimum width=25pt,inner sep=1.5pt]
\tikzstyle{place1}=[rectangle,draw=black,thick,minimum width=20pt,minimum height=20pt]
\node (11) at (275bp,280bp) [place] {$I_{11}$};
  \node (7) at (60bp,105bp) [place] {$I_{7}$};
  \node (2) at (65bp,274bp) [place] {$I_2$};
  \node (3) at (85bp,200bp) [place] {$I_3$};
  \node (1) at (9bp,175bp) [place] {$I_1$};
  \node (5) at (164bp,165bp) [place1] {$I_5$};
  \node (4) at (176bp,244bp) [place] {$I_4$};
  \node (10) at (175bp,85bp) [place] {$I_{10}$};
  \node (6) at (250bp,161bp) [place] {$I_6$};
  \node (9) at (296bp,95bp) [place] {$I_9$};
  \node (8) at (315bp,210bp) [place] {$I_8$};
  \draw [->] (7) -- (5);
  \draw [->] (1) -- (3);
  \draw [->] (10) -- (7);
  \draw [->] (11) -- (8);
  \draw [->] (4) -- (5);
  \draw [->] (3) -- (5);
  \draw [->] (6) -- (5);
  \draw [->] (9) -- (6);
  \draw [->] (2) -- (3);
  \draw [->] (8) -- (6);
\end{tikzpicture}
}
\qquad\qquad
\subfloat[][The $\mathtt{RevLex}$-deformation graph.]{\label{fig:BorelGraph_2}
  \begin{tikzpicture}[>=latex,scale=0.5]
\tikzstyle{place}=[ellipse,draw=black!50,minimum width=25pt,inner sep=1.5pt]
\tikzstyle{place1}=[rectangle,draw=black,thick,minimum width=20pt,minimum height=20pt]
\node (11) at (275bp,280bp) [place1] {$I_{11}$};
  \node (7) at (60bp,105bp) [place] {$I_{7}$};
  \node (2) at (65bp,274bp) [place] {$I_2$};
  \node (3) at (85bp,200bp) [place] {$I_3$};
  \node (1) at (9bp,175bp) [place] {$I_1$};
  \node (5) at (164bp,165bp) [place] {$I_5$};
  \node (4) at (176bp,244bp) [place] {$I_4$};
  \node (10) at (175bp,85bp) [place1] {$I_{10}$};
  \node (6) at (250bp,161bp) [place] {$I_6$};
  \node (9) at (296bp,95bp) [place] {$I_9$};
  \node (8) at (315bp,210bp) [place] {$I_8$};
  \draw [->] (7) -- (10);
  \draw [->] (2) -- (4);
  \draw [->] (4) -- (5);
  \draw [->] (3) -- (5);
  \draw [->] (6) -- (9);
  \draw [->] (9) -- (10);
  \draw [->] (1) -- (3);
  \draw [->] (8) -- (9);
  \draw [->] (5) -- (7);
\end{tikzpicture}
}
\end{center}
\caption{\label{fig:BorelGraph} Two examples of deformation graph of the Hilbert scheme $\hilb^3_{6t-5}$. The square vertices are endpoints.}
\end{figure}
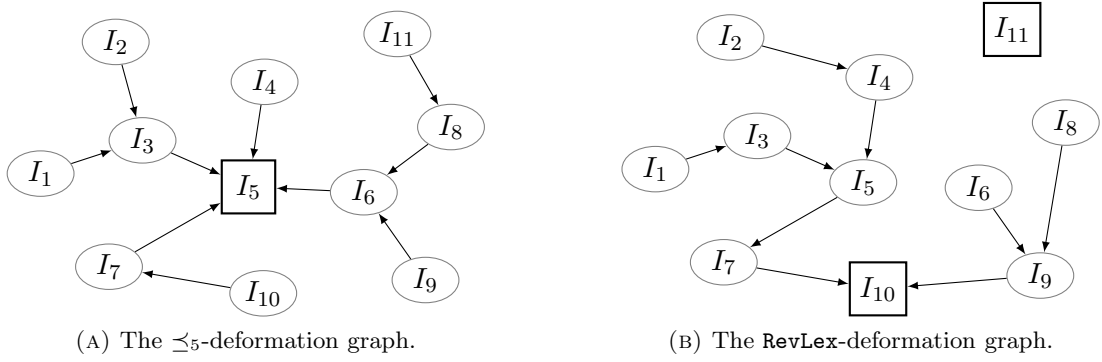 

\end{example}

Let us finally consider the special case of the $\mathtt{DegLex}$-deformation graph of Hilbert schemes of points. Given a Borel-fixed ideal $I \in \hilbd$, we highlight that among the minimal monomials of $\{I_d\}$ there will be surely a monomial of the type $x_1^{s} x_0^{d-s}$. In fact $x_1^d$ belongs to the ideal and applying the decreasing move $\down_1$ repeatedly, we will find a monomial $x_1^{s} x_0^{d-s}$ such that $\down_1(x_1^s x_0^{d-s}) \notin I$. It is easily to deduce that the power $s$ of the variable $x_1$ is the regularity of the saturated ideal $\textnormal{sat}(I)$, because the generators of the saturation of a Borel-fixed ideal can be obtained substituting 1 to the smallest variable $x_0$ and the regularity of a Borel-fixed ideal coincides with the maximal degree of a generator \cite[Proposition 2.9, Proposition 2.11]{Green}.

\begin{proposition}\label{prop:regSteps}
Let $I \in \hilbd$ be a Borel-fixed ideal and let $\textnormal{reg}(I)$ be the regularity of its saturation $\textnormal{sat}(I)$. $d-\textnormal{reg}(I)$ consecutive $\mathtt{DegLex}$-rational deformations lead from $I$ to the lexsegment ideal $L \in \hilbd$.
\begin{proof}
Given any ideal $J \in \hilbd$, we divide the monomials of $N(J)_d$ in two sets:
\[
 \{N(J)_d\} = \{x_0^d,\ldots,x_1^{s-1} x_0^{d-s+1}\} \cup \{x^{\gamma_1},\ldots,x^{\gamma_{d-s}}\},\qquad \max x^{\gamma_i} > 1,\ i=1,\ldots,d-s.
\]
For the lexicographic ideal $L=(x_n,\ldots,x_2,x_1^d)_{\geqslant d}$, the first set contains already $d$ monomials $x_0^d,\ldots,x_1^{d-1} x_0$ and so the second set is empty, whereas the ideal $I$ contains $x_1^{\textnormal{reg}(I)} x_0^{d-\textnormal{reg}(I)}$,and the decomposition of $N(I)_d$ is
\[
 \{N(I)_d\} = \left\{x_0^d,\ldots,x_1^{\textnormal{reg}(I)-1} x_0^{d-\textnormal{reg}(I)+1}\right\} \cup \left\{x^{\gamma_1},\ldots,x^{\gamma_{d-\textnormal{reg}(I)}}\right\}.
\]

Applying Algorithm \ref{alg:deformation} on $I$ and $\mathtt{DegLex}$, we have
\[
 \begin{split}
  \min_{\mathtt{DegLex}} \mathfrak{m}\{I_d\}_0 & = x_1^\delta x_0^{d-\textnormal{reg}(I)},\\
  \max_{\mathtt{DegLex}} \mathfrak{M}\{N(I)_d\}_0 & = \max_{\mathtt{DegLex}} \{x^{\gamma_1},\ldots,x^{\gamma_{d-\delta}}\} = x^{\overline{\gamma}}
 \end{split}
\]
and $x^{\overline{\gamma}} \succ_{\mathtt{DegLex}} x_1^\delta x_0^{d-\textnormal{reg}(I)}$, because for any $x^{\gamma_i}$, $\max x^{\gamma_i} > 1$, and we obtain the ideal ${I}^1 = \left\langle\{I_d\}\setminus\{x_1^\delta x_0^{d-\textnormal{reg}(I)}\}\cup\{x^{\overline{\gamma}}\}\right\rangle$.

$\textnormal{reg}(I^1) = \textnormal{reg}(I)+1$ and repeating $s$ times this process we have $\textnormal{reg}(I^s) = \textnormal{reg}(I) + s$.  
After $d-\textnormal{reg}(I)$ consecutive $\mathtt{DegLex}$-deformations, $\textnormal{reg}(I^{d-\textnormal{reg}(I)}) = d = \textnormal{reg}(L)$, that is $N(I^{d-\textnormal{reg}(I)}) = N(L)$.
\end{proof}
\end{proposition}

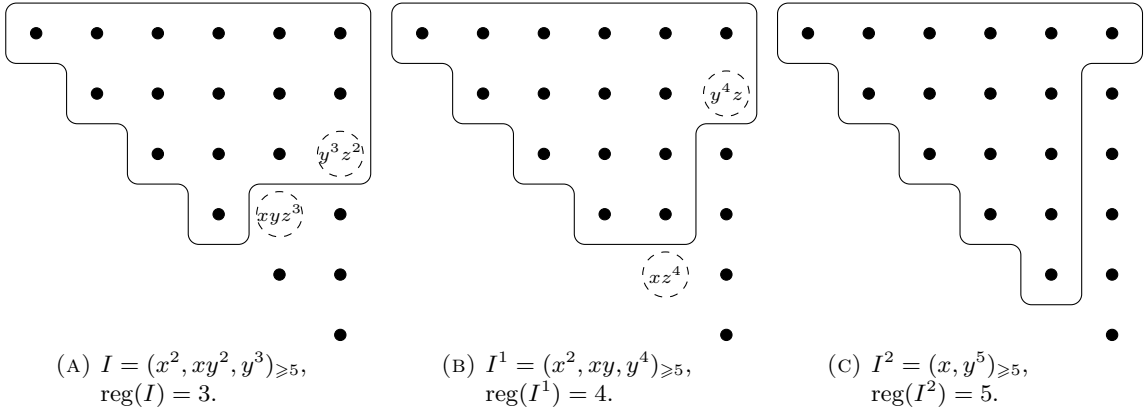
\begin{figure}[!ht]
\begin{center}
\captionsetup[subfloat]{singlelinecheck=false,width=3.5cm,format=hang}
\subfloat[][$I = (x^2,xy^2,y^3)_{\geqslant 5}$,\\ $\textnormal{reg}(I) = 3$.]{
 \begin{tikzpicture}[label distance=-4pt]
\tikzstyle{place}=[circle,draw=black,fill=black,inner sep=1.5pt]
\tikzstyle{place1}=[circle,draw=black,thick,inner sep=1.5pt]

\node (a_1) at (0.4,4.4) [place] {};
\node (a_2) at (1.2,4.4) [place] {};
\node (a_3) at (2,4.4) [place] {};
\node (a_4) at (2.8,4.4) [place] {};
\node (a_5) at (3.6,4.4) [place] {};
\node (a_6) at (4.4,4.4) [place] {};

\node (a_7) at (1.2,3.6) [place] {};
\node (a_8) at (2,3.6) [place] {};
\node (a_9) at (2.8,3.6) [place] {};
\node (a_10) at (3.6,3.6) [place] {};
\node (a_11) at (4.4,3.6) [place] {};

\node (a_12) at (2,2.8) [place] {};
\node (a_13) at (2.8,2.8) [place] {};
\node (a_14) at (3.6,2.8) [place] {};
\node (a_15) at (4.4,2.8) {\tiny $y^3 z^2$};
\draw [dashed] (a_15) circle (0.3);

\node (a_16) at (2.8,2) [place] {};
\node (a_17) at (3.6,2) {\tiny $xyz^3$};
\draw [dashed] (a_17) circle (0.3);
\node (a_18) at (4.4,2) [place] {};

\node (a_19) at (3.6,1.2) [place] {};
\node (a_20) at (4.4,1.2) [place] {};

\node (a_21) at (4.4,0.4) [place] {};

\draw[rounded corners=4pt] (0,4.8) -- (4.8,4.8) -- (4.8,2.4) -- (3.2,2.4) -- (3.2,1.6) -- (2.4,1.6) -- (2.4,2.4) -- (1.6,2.4) -- (1.6,3.2) -- (0.8,3.2) -- (0.8,4) -- (0,4) -- cycle;
\end{tikzpicture}
}
\subfloat[][$I^1 = (x^2,xy,y^4)_{\geqslant 5}$,\\ $\textnormal{reg}(I^1) = 4$.]{
 \begin{tikzpicture}[label distance=-4pt]
\tikzstyle{place}=[circle,draw=black,fill=black,inner sep=1.5pt]
\tikzstyle{place1}=[circle,draw=black,thick,inner sep=1.5pt]

\node (b_1) at (5.6,4.4) [place] {};
\node (b_2) at (6.4,4.4) [place] {};
\node (b_3) at (7.2,4.4) [place] {};
\node (b_4) at (8,4.4) [place] {};
\node (b_5) at (8.8,4.4) [place] {};
\node (b_6) at (9.6,4.4) [place] {};

\node (b_7) at (6.4,3.6) [place] {};
\node (b_8) at (7.2,3.6) [place] {};
\node (b_9) at (8,3.6) [place] {};
\node (b_10) at (8.8,3.6) [place] {};
\node (b_11) at (9.6,3.6) {\tiny $y^4 z$};
\draw [dashed] (b_11) circle (0.3);

\node (b_12) at (7.2,2.8) [place] {};
\node (b_13) at (8,2.8) [place] {};
\node (b_14) at (8.8,2.8) [place] {};
\node (b_15) at (9.6,2.8) [place] {};

\node (b_16) at (8,2) [place] {};
\node (b_17) at (8.8,2) [place] {};
\node (b_18) at (9.6,2) [place] {};

\node (b_19) at (8.8,1.2) {\tiny $xz^4$};
\draw [dashed] (b_19) circle (0.3);
\node (b_20) at (9.6,1.2) [place] {};

\node (b_21) at (9.6,0.4) [place] {};

\draw[rounded corners=4pt] (5.2,4.8) -- (10,4.8) -- (10,3.2) -- (9.2,3.2) -- (9.2,1.6) -- (7.6,1.6) -- (7.6,2.4) -- (6.8,2.4) -- (6.8,3.2) -- (6,3.2) -- (6,4) -- (5.2,4) -- cycle;
\end{tikzpicture}
}
\subfloat[][$I^2 = (x,y^5)_{\geqslant 5}$,\\ $\textnormal{reg}(I^2) = 5$.]{
 \begin{tikzpicture}[label distance=-4pt]
\tikzstyle{place}=[circle,draw=black,fill=black,inner sep=1.5pt]
\tikzstyle{place1}=[circle,draw=black,thick,inner sep=1.5pt]

\node (c_1) at (10.8,4.4) [place] {};
\node (c_2) at (11.6,4.4) [place] {};
\node (c_3) at (12.4,4.4) [place] {};
\node (c_4) at (13.2,4.4) [place] {};
\node (c_5) at (14,4.4) [place] {};
\node (c_6) at (14.8,4.4) [place] {};

\node (c_7) at (11.6,3.6) [place] {};
\node (c_8) at (12.4,3.6) [place] {};
\node (c_9) at (13.2,3.6) [place] {};
\node (c_10) at (14,3.6) [place] {};
\node (c_11) at (14.8,3.6) [place] {};

\node (c_12) at (12.4,2.8) [place] {};
\node (c_13) at (13.2,2.8) [place] {};
\node (c_14) at (14,2.8) [place] {};
\node (c_15) at (14.8,2.8) [place] {};

\node (c_16) at (13.2,2) [place] {};
\node (c_17) at (14,2) [place] {};
\node (c_18) at (14.8,2) [place] {};

\node (c_19) at (14,1.2) [place] {};
\node (c_20) at (14.8,1.2) [place] {};

\node (c_21) at (14.8,0.4) [place] {};

\draw[rounded corners=4pt] (10.4,4.8) -- (15.2,4.8) -- (15.2,4) -- (14.4,4) -- (14.4,0.8) -- (13.6,0.8) -- (13.6,1.6) -- (12.8,1.6) -- (12.8,2.4) -- (12,2.4) -- (12,3.2) -- (11.2,3.2) -- (11.2,4) -- (10.4,4) -- cycle;
\end{tikzpicture}
}
\end{center}
\caption{\label{fig:lexDeformation} An example of Proposition \ref{prop:regSteps}. We represent the two $\mathtt{DegLex}$-rational deformations on $\hilb_5^2$ that lead from the ideal $I = (x^2,xy^2,y^3)_{\geqslant 5}$ to the lexsegment ideal $L = I^2 = (x,y^5)_{\geqslant 5}$.}
\end{figure}

\begin{corollary}\label{cor:treeDepth}
 Given $\hilbd$ the Hilbert scheme of $d$ points in $\PP^n$, let $s$ be the positive integer such that
\[
 \sum_{i=0}^{s-1} \binom{n-1 + i}{n-1} < d \leqslant \sum_{i=0}^{s} \binom{n-1 + i}{n-1}.
\]
The $\mathtt{DegLex}$-deformation graph has height equal to $d-s-1$.
\end{corollary}
We recall that the height of a rooted tree is the maximal distance between a vertex and the root, where the vertices connected to the root have distance 1, the vertices connected to vertices of distance 1 have distance 2 and so on.
\begin{proof}
Let us consider the monomials in $\K[x]_d$ ordered w.r.t. the $\mathtt{RevLex}$ term ordering and let $B$ be the Borel set containing the first $\binom{n+d}{d}-d$ monomials. Obviously the ideal $I = \langle B \rangle$ belongs to $\hilbd$. Moreover by definition of $s$, $x_1^{s} x_0^{d-s}$ is in $N(I)_d$ and $x_1^{s+1} x_0^{d-s-1}$ in in $I_d$, so the regularity of the saturation of $I$ is $s+1$.
\begin{center}
 \begin{tikzpicture}[scale=0.5]
 \draw [rounded corners=3pt,thin,black!50] (-0.1,0.1) -- (-0.1,0.9) -- (-0.9,0.9) -- (-0.9,0.1) -- cycle;
 \draw [rounded corners=3pt,thin,black!50] (-0.1,1.1) -- (-0.1,1.9) -- (-1.9,1.9) -- (-1.9,1.1) -- cycle;
 \draw [rounded corners=3pt,thin,black!50] (-0.1,3.1) -- (-0.1,3.9) -- (-3.9,3.9) -- (-3.9,3.1) -- cycle;
 \draw [rounded corners=3pt,thin,black!50] (-0.1,4.1) -- (-0.1,4.9) -- (-4.9,4.9) -- (-4.9,4.1) -- cycle;
\draw [rounded corners=3pt,thin,black!50] (-0.1,5.1) -- (-0.1,5.9) -- (-5.9,5.9) -- (-5.9,5.1) -- cycle;
\draw [rounded corners=3pt,thin,black!50] (-0.1,7.6) -- (-0.1,8.4) -- (-8.4,8.4) -- (-8.4,7.6) -- cycle;
\draw [thick] (0,5) -- (-3,5) -- (-3,4) -- (-5,4);
\node at (-0.5,2.65) [] {$\vdots$};
\node at (-0.5,6.9) [] {$\vdots$};

\node (e) at (1.1,0.5) [] {\small $x_0^d$};
\node at (3.8,1.5) {\small $\K[x_1,\ldots,x_n]_1\cdot x_0^{d-1}$};
\node at (4.45,3.5) [] {\small $\K[x_1,\ldots,x_n]_{s+1}\cdot x_0^{d-s-1}$};
\node at (3.7,4.5) [] {\small $\K[x_1,\ldots,x_n]_{s}\cdot x_0^{d-s}$};
\node at (4.45,5.5) [] {\small $\K[x_1,\ldots,x_n]_{s-1}\cdot x_0^{d-s+1}$};
\node at (2.8,8) [] {\small $\K[x_1,\ldots,x_n]_d$};

\node at (-8,6.5) [] {Monomials in $I_d$};
\node at (-6,2.4) [] {Monomials in $N(I)_d$};
 \end{tikzpicture}
\end{center}

By Proposition \ref{prop:regSteps}, we know that from $I$ we can reach the lexicographic ideal $L$ (root of the $\mathtt{DegLex}$-deformation graph) through $d-s-1$ rational deformation. Finally by definition of $\mathtt{RevLex}$, for any Borel ideal $J \in \hilbd$, $\textnormal{reg}(J) \geqslant \textnormal{reg}(I)$, so there does not exist an ideal more distant from $L$ than $I$.
\end{proof}

Corollary \ref{cor:treeDepth} seems to suggest that hilb-segment ideals, such that their saturation has high or low regularity, are in a certain sense in the \lq\lq suburb\rq\rq of the Hilbert scheme $\hilbd$. In fact called $\delta$ the regularity of the saturation of the hilb-segment w.r.t. $\mathtt{RevLex}$, the corresponding deformation graphs should have a height in the vicinity of $d - \delta$. Whereas the hilb-segment ideals with saturation having regularity near to $\frac{d+\delta}{2}$ seems to be in the \lq\lq heart\rq\rq\ of the Hilbert scheme, because the corresponding deformation graph should have minimal height.

\begin{example}\label{ex:depthDefGraph}
 Let us consider the Hilbert scheme $\hilb_8^3$. There are 12 Borel-fixed ideals:
\[
 \begin{array}{l cc l}
  J_1 = (x_3,x_2,x_1^8)_{\geqslant 8}, &&& J_7 = (x_3^2,x_3x_2,x_2^2,x_3x_1,x_2x_1^2,x_1^5)_{\geqslant 8}, \\
  J_2 = (x_3,x_2^2,x_2x_1,x_1^7)_{\geqslant 8} &&& J_8 = (x_3,x_2^3,x_2^2x_1,x_2x_1^3,x_1^4)_{\geqslant 8}, \\
  J_3 = (x_3,x_2^2,x_2x_1^2,x_1^6)_{\geqslant 8}, &&& J_9 = (x_3^2,x_3x_2,x_2^2,x_3x_1,x_2x_1^3,x_1^4)_{\geqslant 8},\\
  J_4 = (x_3^2,x_3x_2,x_2^2,x_3x_1,x_2x_1,x_1^6)_{\geqslant 8}, &&&  J_{10} = (x_3^2,x_3x_2,x_3x_1,x_2^3,x_2^2x_1,x_2x_1^2,x_1^4)_{\geqslant 8},\\
  J_5 = (x_3,x_2^2,x_2x_1^3,x_1^5)_{\geqslant 8}, &&& J_{11} = (x_3^2,x_3x_2,x_2^2,x_3x_1^2,x_2x_1^2,x_1^4)_{\geqslant 8},\\
  J_6 = (x_3,x_2^3,x_2^2x_1,x_2x_1^2,x_1^5)_{\geqslant 8}, &&& J_{12} = (x_3^2,x_3x_2,x_2^3,x_2^2x_1,x_3x_1^2,x_2x_1^2,x_1^3)_{\geqslant 8}.
 \end{array}
\]
Since $\binom{2}{2} + \binom{2+1}{2} < 8 < \binom{2}{2} + \binom{2+1}{2} + \binom{2+2}{2}$, by Corollary \ref{cor:treeDepth} we know that the $\mathtt{DegLex}$-deformation graph has height $5$. $5$ is also the height of the $\mathtt{RevLex}$-deformation graph (Figure \ref{fig:depthDefGraph}\subref{fig:depthDefGraph_1}), that has $J_{12}$ as root. 

$J_7$ is a hilb-segment w.r.t. the term ordering defined by the matrix in \eqref{eq:matrixTO} with $\omega=(5,4,2,1)$ and its saturation has regularity $5$. The $\preceq_7$-deformation graph has minimal height equal to $3$ (Figure \ref{fig:depthDefGraph}\subref{fig:depthDefGraph_2}).
\begin{figure}[!ht]
 \begin{center}\captionsetup[subfloat]{width=5cm}
\subfloat[][The $\mathtt{RevLex}$-deformation graph.]{\label{fig:depthDefGraph_1}
\begin{tikzpicture}[>=latex,scale=0.85]
\tikzstyle{place}=[ellipse,draw=black!50,minimum width=25pt,inner sep=1.5pt]
\tikzstyle{place1}=[rectangle,draw=black,thick,minimum width=20pt,minimum height=20pt]
  \node (11) at (74bp,160bp) [place] {$J_{11}$};
  \node (10) at (116bp,160bp) [place] {$J_{10}$};
  \node (12) at (95bp,200bp) [place1] {$J_{12}$};
  \node (1) at (8bp,0bp) [place] {$J_1$};
  \node (3) at (42bp,80bp) [place] {$J_3$};
  \node (2) at (8bp,40bp) [place] {$J_2$};
  \node (5) at (76bp,80bp) [place] {$J_5$};
  \node (4) at (8bp,80bp) [place] {$J_4$};
  \node (7) at (42bp,120bp) [place] {$J_7$};
  \node (6) at (114bp,120bp) [place] {$J_6$};
  \node (9) at (76bp,120bp) [place] {$J_9$};
  \node (8) at (148bp,120bp) [place] {$J_8$};
  \draw [->] (3) -- (7);
  \draw [->] (5) -- (9);
  \draw [->] (4) -- (7);
  \draw [->] (10) -- (12);
  \draw [->] (8) -- (10);
  \draw [->] (6) -- (10);
  \draw [->] (1) -- (2);
  \draw [->] (9) -- (11);
  \draw [->] (11) -- (12);
  \draw [->] (7) -- (11);
  \draw [->] (2) -- (4);
\end{tikzpicture}
}
\quad\quad
\subfloat[][The $\preceq_7$-deformation graph.]{\label{fig:depthDefGraph_2}
\begin{tikzpicture}[>=latex,scale=0.85]
\tikzstyle{place}=[ellipse,draw=black!50,minimum width=25pt,inner sep=1.5pt]
\tikzstyle{place1}=[rectangle,draw=black,thick,minimum width=20pt,minimum height=20pt]
  \node (11) at (224bp,160bp) [place] {$J_{11}$};
  \node (10) at (182bp,160bp) [place] {$J_{10}$};
  \node (12) at (224bp,120bp) [place] {$J_{12}$};
  \node (1) at (8bp,80bp) [place] {$J_1$};
  \node (3) at (42bp,160bp) [place] {$J_3$};
  \node (2) at (8bp,120bp) [place] {$J_2$};
  \node (5) at (76bp,160bp) [place] {$J_5$};
  \node (4) at (8bp,160bp) [place] {$J_4$};
  \node (7) at (110bp,200bp) [place1] {$J_7$};
  \node (6) at (110bp,160bp) [place] {$J_6$};
  \node (9) at (144bp,160bp) [place] {$J_9$};
  \node (8) at (144bp,120bp) [place] {$J_8$};
 \node at (0bp,0bp) [draw=white] {$\phantom{J_8}$};
  \draw [->] (11) -- (7);
  \draw [->] (3) -- (7);
  \draw [->] (10) -- (7);
  \draw [->] (5) -- (7);
  \draw [->] (4) -- (7);
  \draw [->] (6) -- (7);
  \draw [->] (12) -- (11);
  \draw [->] (9) -- (7);
  \draw [->] (1) -- (2);
  \draw [->] (8) -- (9);
  \draw [->] (2) -- (4);
\end{tikzpicture}
}
 \end{center}
 \caption{\label{fig:depthDefGraph} The deformation graphs of Example \ref{ex:depthDefGraph} viewed as rooted trees.}
\end{figure}
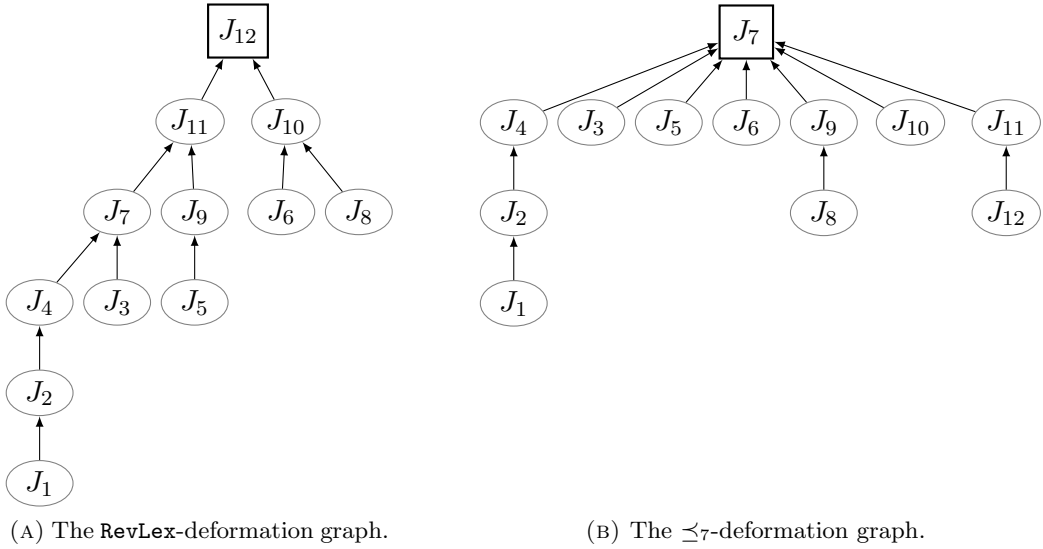
\end{example}

\section{Components of the Hilbert scheme}\label{sec:components}
\begin{figure}[!ht] 
\begin{center}
\subfloat[][Example \ref{ex:PSsameRational}.]{
\begin{tikzpicture}[>=latex,scale=1.05]
\tikzstyle{place}=[circle,draw,fill,inner sep=0pt,minimum size=1.5mm]
\tikzstyle{place1}=[shape=ellipse,draw,minimum width=5.2cm,minimum height=1.8cm]
\node (I) at (0,0) [place,label=270:{\small $I$}] {};
\node (I2) at (2,3) [place,label=90:{\small $\widetilde{I}$}] {};
\node (J) at (4,0) [place,label=270:{\small $J$}] {};

\draw [->,thin,decorate,decoration={snake,amplitude=.3mm,segment length=2mm}] (I2) --node[fill=white]{\tiny $\text{in}_{\mathtt{RevLex}}$} (I);
\draw [->,thin,decorate,decoration={snake,amplitude=.3mm,segment length=2mm}] (I2) --node[fill=white]{\tiny $\text{in}_{\mathtt{Lex}}$} (J);
\draw [<->,thin] (I) --node[below]{\tiny $\mathtt{DegLex}$-deformation} (J);

\node at (3,1.5) [place1,rotate=124,very thin,black!50] {};
\node at (1,1.5) [place1,rotate=56,very thin,black!50] {};
\draw [thin,dashed,black!50] (2,0) ellipse (2.5 and 0.8);
\end{tikzpicture}
}
\qquad\qquad
\subfloat[][Example \ref{ex:PSbreakingBorel}.]{
  \begin{tikzpicture}[>=latex]
\tikzstyle{place}=[circle,draw,fill,inner sep=0pt,minimum size=1.5mm]
\tikzstyle{place1}=[shape=ellipse,draw,minimum width=4.8cm,minimum height=1.6cm]

\node (I) at (0,0) [place,label=270:{\small $I$}] {};
\node (I2) at (1,3) [place,label=90:{\small $\widetilde{I}$}] {};
\node (I3) at (3,3) [place,label=90:{\small $\text{in}_{\mathtt{Lex}}(\widetilde{I})$}] {};
\node (J) at (4,0) [place,label=270:{\small $J$}] {};

\draw [<->,thin] (I) --node[below]{\tiny $\mathtt{DegLex}$-deformation} (J);
\draw [->,thin,decorate,decoration={snake,amplitude=.3mm,segment length=2mm}] (I2) --node[fill=white]{\tiny $\text{in}_{\mathtt{RevLex}}$} (I);
\draw [->,thin,decorate,decoration={snake,amplitude=.3mm,segment length=2mm}] (I2) --node[fill=white]{\tiny $\text{in}_{\mathtt{Lex}}$} (I3);
\draw [->,thin,decorate,decoration={snake,amplitude=.3mm,segment length=2mm}] (I3) --node[fill=white]{\tiny $\text{gin}_{\mathtt{Lex}}$} (J);

\draw [thin,dashed,black!50] (2,0) ellipse (2.5 and 0.8);
\draw [very thin,black!50] (2,3) ellipse (1.7 and 0.7);
\node at (3.5,1.5) [place1,rotate=108,very thin,black!50] {};
\node at (0.5,1.5) [place1,rotate=72,very thin,black!50] {};

\end{tikzpicture}
}
\end{center}
\caption{\label{fig:ComponentGraph} Applying our rational deformation, we can know for sure whenever two Borel ideals lie on a same component (drawn with the dashed line). With Peeva-Stillman method, in the best case (on the left), we know that the components containing two Borel ideals have non-empty intersection but in the other cases, we could get no information (on the right).}
\end{figure}
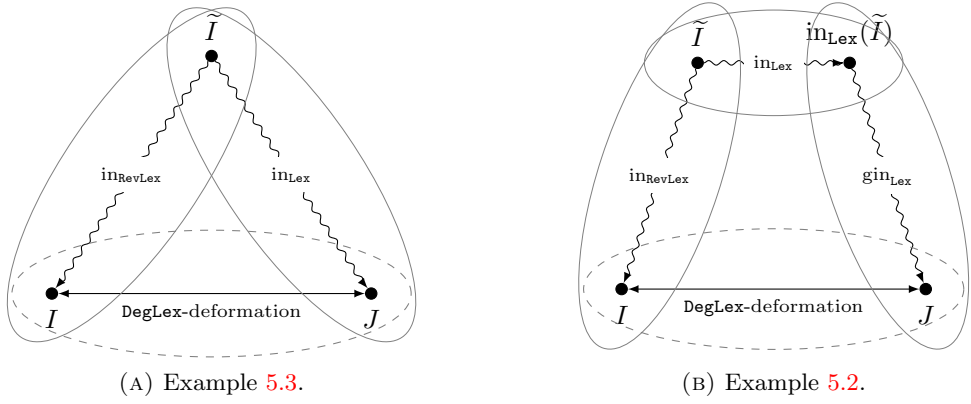

The last section of the paper is devoted to the study of components of the Hilbert scheme, keeping in mind a question posed by Reeves in \cite{RA} \lq\lq Is the subset of Borel-fixed ideals on a component enough to determine the component?\rq\rq. Using our method in fact the points corresponding to two ideals connected by any rational deformation lie on the same components. We underline that the technique introduced in \cite{PS} is slightly different, because in order to pass from a Borel ideal to another Peeva and Stillman use at least two affine deformation (see Figure \ref{fig:ComponentGraph}). 

We introduce a new graph related to Borel-fixed ideals and Hilbert scheme components, different from the incidence graph defined by Reeves \cite{RA}, but also useful to understand the intersections among components and we will call it Borel incidence graph.
\begin{definition}
 Given the Hilbert scheme $\hilbp$, we define the \emph{Borel incidence graph} as the graph $(V,E)$ such that:
\begin{itemize}
 \item the set of vertices $V$ contains the Borel-fixed ideals defining points of $\hilbp$;
 \item an edge $(I,J)$ belongs to $E$ whenever a rational deformation connecting $I$ and $J$ exists.
\end{itemize}
\end{definition}

To construct it we can apply Algorithm \ref{alg:allDeformations} on each Borel ideal of $\hilbp$ and put together all the rational deformations found (discarding repetitions). But we want to add more edges to the graph, so we now discuss if two rational deformations can be performed simultaneously. Let us consider an ideal $I \in \hilbp$ having two rational deformations: $\phi_1$ leading to the ideal $J^1 = \left\langle\{I\}\setminus\mathcal{F}_{\alpha_1}(x^{\alpha_1})\cup\mathcal{F}_{\alpha_1}(x^{\beta_1})\right\rangle$ and $\phi_2$ leading to $J^2 = \big\langle\{I\}\setminus\mathcal{F}_{\alpha_2}(x^{\alpha_2})$ $\left.\cup\mathcal{F}_{\alpha_2}(x^{\beta_2})\right\rangle$. The point is to understand whenever swapping at the same time $\mathcal{F}_{\alpha_1}(x^{\alpha_1})$ with $\mathcal{F}_{\alpha_1}(x^{\beta_1})$ and $\mathcal{F}_{\alpha_2}(x^{\alpha_2})$ with $\mathcal{F}_{\alpha_2}(x^{\beta_2})$, the Borel condition can be broken.
\begin{enumerate}
 \item First of all if the sets of monomials are not disjointed, there are some problem in the definition of the deformation. For instance if $\mathcal{F}_{\alpha_1}(x^{\alpha_1}) \cap \mathcal{F}_{\alpha_2}(x^{\alpha_2}) = \emptyset$ and $\mathcal{F}_{\alpha_1}(x^{\beta_1}) \cap \mathcal{F}_{\alpha_2}(x^{\beta_2}) \neq \emptyset$, the number of monomials to move from the ideal to the quotient would be different from the monomials moving in the opposite direction.
 \item If $x^{\beta_2}$ can be obtained by decreasing moves from $x^{\alpha_1}$ (or $x^{\beta_1}$ from $x^{\alpha_2}$), performing both swaps we would obtain an ideal $\widetilde{I}$ not Borel, because with an element in it Borel-smaller than a monomial in its quotient: $\widetilde{I} \ni x^{\beta_2} \leq_B x^{\alpha_1} \notin \widetilde{I}$.
 \item Let us suppose $\min x^{\alpha_1} > \min x^{\alpha_2}$. If there exists an admissible move $\down_k$, $k \geqslant \min x^{\alpha_1}$, such that $\down_k(x^\alpha_2) \leq_B x^{\beta_1}$, it could happen that $\down_k(x^{\alpha_2}) \in \mathcal{F}_{\alpha_1}(x^{\beta_1})$ and again swapping both sets we would obtain an ideal $\widehat{I}$ not Borel because $\widehat{I} \ni \down_k(x^{\alpha_2}) \leq_B x^{\alpha_2} \notin \widehat{I}$.
\end{enumerate}

\begin{theorem}\label{th:compatible}
 Let $I \in \hilbp$ be a Borel-fixed ideal having $s$ rational deformations: 
\[
\phi_k: I \rightsquigarrow J_k = \left\langle\{I\}\setminus\mathcal{F}_{\alpha_k}(x^{\alpha_k})\cup\mathcal{F}_{\alpha_k}(x^{\beta_k})\right\rangle,\ x^{\alpha_k} \in \mathfrak{m}\{I\}_{h_k},\ x^{\beta_k} \in \mathfrak{M}\{\NI\}_{h_k},
\]
$k=1,\ldots,s$, $\min x^{\alpha_1} \geqslant \ldots \geqslant \min^{\alpha_s}$.
If the following conditions hold:
\begin{enumerate}[(i)]
 \item\label{it:compatible_hp1} $\forall\ i,j$, $x^{\alpha_i}$ can not be obtained from $x^{\beta_j}$ by applying a decreasing move;
 \item\label{it:compatible_hp2} $\forall\ i,j$ such that $\min x^{\alpha_i} > \min x^{\alpha_j}$, $\down_l(x^{\alpha_j})$ cannot be obtained from $x^{\beta_i}$ through a decreasing move, for all admissible $\down_l,\ l \geqslant \min x^{\alpha_i}$;
 \item\label{it:compatible_hp3} $\forall\ i,j$, $\mathcal{F}_{\alpha_i}(x^{\alpha_i}) \cap \mathcal{F}_{\alpha_j}(x^{\alpha_j}) = \emptyset$ and $\mathcal{F}_{\alpha_i}(x^{\beta_i}) \cap \mathcal{F}_{\alpha_j}(x^{\beta_j}) = \emptyset$;
\end{enumerate}
then there exists a flat family over $\left(\PP^1\right)^{\times s}$ of subschemes parametrized by $\hilbp$ containing $2^s$ Borel-fixed ideals.
\begin{proof}
It suffices to repeat the reasoning of the proof of Lemma \ref{lem:BorderMonomials} and Theorem \ref{th:flatGeneralDim}, starting from the family of ideals
\begin{equation}\label{eq:multipleDef}
\begin{split}
 \mathcal{I} = \Big\langle \left\{I\right\} &\setminus \mathcal{F}_{\alpha_1}(x^{\alpha_1})\cup\left\{y_{10}\, F(x^{\alpha_1}) + y_{11}\, F(x^{\beta_1})\ \big\vert\ F \in\mathcal{F}_{\alpha_1},\ [y_{10}:y_{11}] \in \PP^1\right\}\\
  & \ldots\\
  &\setminus \mathcal{F}_{\alpha_i}(x^{\alpha_i})\cup\left\{y_{i0}\, F(x^{\alpha_i}) + y_{i1}\, F(x^{\beta_i})\ \big\vert\ F \in\mathcal{F}_{\alpha_i},\ [y_{i0}:y_{i1}] \in \PP^1\right\}\\
  & \ldots\\
  &\setminus \mathcal{F}_{\alpha_s}(x^{\alpha_s})\cup\left\{y_{s0}\, F(x^{\alpha_s}) + y_{s1}\, F(x^{\beta_s})\ \big\vert\ F \in\mathcal{F}_{\alpha_s},\ [y_{s0}:y_{s1}] \in \PP^1\right\}\Big\rangle.
\end{split}
\end{equation}
For instance in the case $y_{10} \neq 0,\ldots, y_{s0} \neq 0$, a basis of the ideal in degree $r+1$ will be given by the basis $\mathcal{G}'$ of the monomials of degree $r+1$ belonging to the ideal generated by $\left\langle\{I\}\setminus\left(\bigcup\limits_{i=1}^s\mathcal{F}_{\alpha_i}(x^{\alpha_i})\right)\right\rangle$ and the set of monomials
\[
 \mathcal{G}'' = \bigcup_{i=1}^s \left\{ x_j\, F(x^{\alpha_i}) + y_{i1}\, x_j\, F(x^{\beta_i})\ \big\vert\ F\in\mathcal{F}_{\alpha_i}\text{ and }j=0,\ldots,\min F(x^{\alpha_i})\right\}.
\]
The family is given by the morphism of rings 
\[
\widehat{\phi}: \K[y_{10},y_{11}]\otimes_{\K}\ldots\otimes_{\K} \K[y_{s0},y_{s1}] \rightarrow \K[x][y_{10},y_{11},\ldots,y_{s0},y_{s1}]/\mathcal{I}.
\]
Choosing in any way for every $i=1,\ldots,s$ one set of monomials between $\mathcal{F}_{\alpha_i}(x^{\alpha_i})$ and $\mathcal{F}_{\alpha_i}(x^{\beta_i})$, we obtain a Borel-fixed ideal belonging to the family $\mathcal{I}$. In fact such Borel ideals correspond to the points $(\ldots,[y_{i0}:y_{i1}],\ldots) \in \left(\PP^1\right)^{\times s}$, $[y_{i0}:y_{i1}] = [0:1]$ or $[y_{i0}:y_{i1}] = [1:0]$, and so they are $2^s$.
\end{proof}
\end{theorem}

\begin{corollary}\label{cor:curveCompatible}
 Let $I \in \hilbp$ be a Borel-fixed ideal and let $\phi_1,\ldots,\phi_s$ be $s$ rational deformations of $I$ as in Theorem \ref{th:compatible}. For each couple of Borel ideals $J,\widetilde{J}$ belonging to the family over $\left(\PP^1\right)^{\times s}$, there exists a rational curve connecting them.
\begin{proof}
 Let $\mathcal{I}$ be the family defined in \eqref{eq:multipleDef} and $\mathcal{F}_{\alpha_i}(x^{\alpha_i})$,$\mathcal{F}_{\alpha_i}(x^{\beta_i})$ be the monomials involved in the deformation $\phi_i$, $i=1,\ldots,s$.
 We consider firstly the morphism of ring
 \[
  \K[y_{10},y_{11}]\otimes_{\K}\ldots\otimes_{\K} \K[y_{s0},y_{s1}] \rightarrow \K[y_{10},y_{11}]/C_1 \otimes_{\K}\ldots\otimes_{\K} \K[y_{s0},y_{s1}]/C_s
 \]
where $C_i$ is the ideal
\begin{itemize}
 \item $(0)$, if $\mathcal{F}_{\alpha_i}(x^{\alpha_i})$ belong to $J$ and $\mathcal{F}_{\alpha_i}(x^{\beta_i})$ to $\widetilde{J}$ or vice versa if $\mathcal{F}_{\alpha_i}(x^{\alpha_i})$ belong to $\widetilde{J}$ and $\mathcal{F}_{\alpha_i}(x^{\beta_i})$ to $J$;
 \item $(y_{i1})$, if $\mathcal{F}_{\alpha_i}(x^{\alpha_i})$ belongs to both ideals;
 \item $(y_{i0})$, if $\mathcal{F}_{\alpha_i}(x^{\beta_i})$ belongs to both ideals.
\end{itemize}
In this way we eliminate, the copies of $\PP^1$ corresponding to couples of monomials that do not have to be swapped:
\[
\begin{split}
 &\K[y_{10},y_{11}]/C_1 \otimes_{\K}\ldots\otimes_{\K} \K[y_{s0},y_{s1}]/C_s =\\
 &= \left(\bigotimes_{C_i \neq (0)} \K[y_{i0},y_{i1}]/C_i \right) \otimes_\K \left(\bigotimes_{C_j = (0)} \K[y_{j0},y_{j1}]/C_j \right) \simeq \bigotimes_{C_j = (0)} \K[y_{j0},y_{j1}].
\end{split}
\]
Finally we equals the variables corresponding to monomials moving outside the ideal and those corresponding to monomials moving inside, that is 
\[
 \bigotimes_{C_j = (0)} \K[y_{j0},y_{j1}] \rightarrow \K[\ldots,y_{j0},y_{j1},\ldots]/C \simeq \K[y_0,y_1]
\] 
where the ideal $C$ is generated by the linear forms:
\begin{itemize}
 \item $y_{i0}-y_{k0},y_{i1}-y_{k1}$, if $\mathcal{F}_{\alpha_i}(x^{\alpha_i}),\mathcal{F}_{\alpha_k}(x^{\alpha_k})$ belong to $J$ and $\mathcal{F}_{\alpha_i}(x^{\beta_i}),\mathcal{F}_{\alpha_k}(x^{\beta_k})$ to $\widetilde{J}$;
 \item $y_{i0}-y_{k1},y_{i1}-y_{k0}$, if $\mathcal{F}_{\alpha_i}(x^{\alpha_i}),\mathcal{F}_{\alpha_k}(x^{\beta_k})$ belong to $J$ and $\mathcal{F}_{\alpha_i}(x^{\beta_i}),\mathcal{F}_{\alpha_k}(x^{\alpha_k})$ to $\widetilde{J}$.
\end{itemize}
Composing the two morphism of rings, we obtain a rational curve $\mathcal{C}: \PP^1 \rightarrow (\PP^1)^{\times s}$, passing through the points defined by $J$ and $\widetilde{J}$.
\end{proof}
\end{corollary}

\begin{definition}
 Let $I$ be a Borel ideal in $\hilbp$ and let $\phi_1,\ldots,\phi_s$ be a set of rational deformation of $I$. If the monomials involved in the deformations respect the hypothesis of Theorem \ref{th:compatible} we say that $\{\phi_1,\ldots,\phi_s\}$ is a set of \emph{compatible} deformations.
\end{definition}

Theorem \ref{th:compatible} and Corollary \ref{cor:curveCompatible} introduce many other rational deformations between Borel-fixed ideals. Hence to compute the Borel incidence graph (Algorithm \ref{alg:BorelIncidenceGraph}) we have to add these composed rational deformations to those obtained applying Algorithm \ref{alg:allDeformations} on every ideal.

\begin{algorithm}
 \caption{\label{alg:BorelIncidenceGraph} How to compute the Borel incidence graph. A trial version of this algorithm is available at the web page\\ \href{http://www.dm.unito.it/dottorato/dottorandi/lella/HSC/borelIncidenceGraph.html}{www.dm.unito.it/dottorato/dottorandi/lella/HSC/borelIncidenceGraph.html}}
 \begin{algorithmic}
 \Require $\hilbp$, Hilbert scheme parametrizing subschemes of $\PP^n$ with Hilbert polynomial $p(t)$.
 \Ensure the Borel incidence graph $(V,E)$.
 \Procedure{BorelIncidenceGraph}{$\hilbp$}
\State $V\, \leftarrow$ Borel-fixed ideals in $\K[x]$ with Hilbert polynomial $p(t)$.
\State $E\, \leftarrow \emptyset$;
\ForAll{$I \in V$}
\State $S \leftarrow$ \Call{RationalDeformations}{$I$};
   \ForAll{$J \in S$}
     \If{$(I,J) \notin E$}
     \State $E \leftarrow E \cup \{(I,J)\}$;
     \EndIf
   \EndFor
   \ForAll{$\{\phi_1,\ldots,\phi_s\} \subseteq S,\ s > 1$} \Comment{Looking for compatible rational deformations.}
   \If{$\phi_1,\ldots,\phi_s$ are compatible}
      \If{$(I,(\phi_s \circ \cdots \circ \phi_1)(I)) \notin E$}
         \State $E \leftarrow E \cup \{(I,(\phi_s \circ \cdots \circ \phi_1)(I))\}$;
      \EndIf
   \EndIf
   \EndFor
\EndFor
\State\Return $(V,E)$;
 \EndProcedure
 \end{algorithmic}
\end{algorithm}

\begin{example}
The two deformations
\[
 \begin{split}
  \phi_1 &: I \rightsquigarrow J_1 = \left\langle\{I\}\setminus\{x_3^2 x_0^6\}\cup\{x_2^3 x_0^5\}\right\rangle \\
  \phi_2 &: I \rightsquigarrow J_5 = \left\langle\{I\}\setminus\mathcal{F}(x_2^2 x_1^6)\cup\mathcal{F}(x_3 x_1^7)\right\rangle
 \end{split}
\]
described in Example \ref{ex:multipleDef} are in the hypothesis of Theorem \ref{th:compatible}, so there is a flat family over $\PP^1 \times \PP^1$ defined by the ideal
\[
\begin{split}
\mathcal{I} = \Big\langle& \{I\} \setminus \left\{x_3^2 x_0^6\right\} \cup \left\{ y_0\, x_3^2 x_0^6 + y_1\, x_2^3 x_0^5\right\} \\
& \setminus \mathcal{F}(x_2^2 x_1^6) \cup\left\{z_0\, F(x_2^2 x_1^6) + z_1\, F(x_3 x_1^7)\ \vert\ F \in \mathcal{F}\right\}\Big\rangle,\quad ([y_0:y_1],[z_0:z_1]) \in \PP^1 \times \PP^1.
\end{split}
\]
The Borel ideals belonging to the family are $I,J_1,J_5$ and $\widetilde{I} = (x_3^3,x_3^2 x_2,x_3 x_2^2,x_2^3,x_3^2 x_1, x_3 x_2 x_1,$ $x_3 x_1^3)_{\geqslant 8}$. Let us apply Corollary \ref{cor:curveCompatible} to determine the rational deformations between couples of ideals:
\[
 \begin{split}
  &I \rightsquigarrow J_1,\quad \widehat{\mathcal{C}}_1: \K[y_0,y_1]\otimes\K[z_0,z_1] \rightarrow \K[y_0,y_1]\otimes\K[z_0,z_1]/(z_1) \simeq \K[y_0,y_1], \\
  &I \rightsquigarrow J_5,\quad \widehat{\mathcal{C}}_2: \K[y_0,y_1]\otimes\K[z_0,z_1] \rightarrow \K[y_0,y_1]/(y_1)\otimes\K[z_0,z_1] \simeq \K[z_0,z_1], \\
  &J_1 \rightsquigarrow \widetilde{I},\quad \widehat{\mathcal{C}}_3: \K[y_0,y_1]\otimes\K[z_0,z_1] \rightarrow \K[y_0,y_1]/(y_0)\otimes\K[z_0,z_1] \simeq \K[z_0,z_1], \\
  &J_5 \rightsquigarrow \widetilde{I},\quad \widehat{\mathcal{C}}_4: \K[y_0,y_1]\otimes\K[z_0,z_1] \rightarrow \K[y_0,y_1]\otimes\K[z_0,z_1]/(z_0) \simeq \K[y_0,y_1], \\
  &I \rightsquigarrow \widetilde{I},\quad \widehat{\mathcal{C}}_5:\K[y_0,y_1]\otimes\K[z_0,z_1] \rightarrow \K[y_0,y_1,z_0,z_1]/(y_0-z_0,y_1-z_1) \simeq \K[y_0,y_1],\\
  &J_1 \rightsquigarrow J_5,\quad \widehat{\mathcal{C}}_6:\K[y_0,y_1]\otimes\K[z_0,z_1] \rightarrow \K[y_0,y_1,z_0,z_1]/(y_0-z_1,y_1-z_0) \simeq \K[y_0,y_1].
 \end{split}
\]
The first four deformations are simple rational deformation that can be found using Algorithm \ref{alg:allDeformations}, whereas the fifth and the sixth one are composed.
\end{example}

Theorem \ref{th:compatible} gives us a new criterion to detect points defined by Borel ideals lying on a common component of the Hilbert scheme as the following example shows.

\begin{example}
Let us consider the Hilbert scheme containing the rational normal curve of degree 4, that is $\hilb_{4t+1}^4$. On it there are $12$ Borel-fixed ideals:
\begin{enumerate}
\item[] $J_1 = (x_4,x_3,x_2^5,x_2^4x_1^3)_{\geqslant 7}$,
\item[] $J_2 = (x_4,x_3,x_2^6,x_2^5x_1,x_2^4x_1^2)_{\geqslant 7}$,
\item[] $J_3 = (x_4,x_3^2,x_3x_2,x_3x_1,x_2^5,x_2^4x_1^2)_{\geqslant 7}$,
\item[] $J_4 = (x_4,x_3^2,x_3x_2,x_3x_1^2,x_2^5,x_2^4x_1)_{\geqslant 7}$,
\item[] $J_5 = (x_4^2,x_4x_3,x_3^2,x_4x_2,x_3x_2,x_4x_1,x_3x_1,x_2^5,x_2^4x_1)_{\geqslant 7}$,
\item[] $J_6 = (x_4,x_3^2,x_3x_2,x_2^4,x_3x_1^3)_{\geqslant 7}$,
\item[] $J_7 = (x_4,x_3^2,x_3x_2^2,x_3x_2x_1,x_3x_1^2,x_2^4)_{\geqslant 7}$,
\item[] $J_8 = (x_4^2,x_4x_3,x_3^2,x_4x_2,x_3x_2,x_4x_1,x_3x_1^2,x_2^4)_{\geqslant 7}$,
\item[] $J_9 = (x_4,x_3^2,x_3x_2,x_2^4,x_2^3x_1)_{\geqslant 7}$,
\item[] $J_{10} = (x_4,x_3^2,x_3x_2^2,x_2^3,x_3x_2x_1)_{\geqslant 7}$,
\item[] $J_{11} = (x_4^2,x_4x_3,x_3^2,x_4x_2,x_3x_2,x_4x_1,x_2^3)_{\geqslant 7}$,
\item[] $J_{12} = (x_4^2,x_4x_3,x_3^2,x_4x_2,x_3x_2,x_2^2)_{\geqslant 7}$.
\end{enumerate}
For a result by Reeves \cite{RA}, it is well known that Borel ideals with the same hyperplane section lie on a common component of the Hilbert scheme, that is the ideals from $J_1$ to $J_{8}$ lie on a component and the same hold for $J_{9},J_{10},J_{11}$ and for $J_{12}$. Computing with Algorithm \ref{alg:BorelIncidenceGraph} the Borel incidence graph of $\hilb^4_{4t+1}$ (Figure \ref{fig:BorelIncidence}), we find one family over $\PP^1 \times \PP^1$, containing the ideals $J_7,J_8,J_{10},J_{11}$, so that they lie on a same component.
\begin{figure}[!ht]
\begin{center}
 \begin{tikzpicture}[>=latex,scale=0.6]
\tikzstyle{place}=[ellipse,draw=black!50,minimum width=25pt,inner sep=1.5pt]
\node (11) at (111bp,85bp) [place] {$J_{11}$};
  \node (10) at (211bp,83bp) [place] {$J_{10}$};
  \node (12) at (20bp,30bp) [place] {$J_{12}$};
  \node (1) at (431bp,332bp) [place] {$J_1$};
  \node (3) at (409bp,256bp) [place] {$J_3$};
  \node (2) at (552bp,292bp) [place] {$J_2$};
  \node (5) at (414bp,179bp) [place] {$J_5$};
  \node (4) at (277bp,211bp) [place] {$J_4$};
  \node (7) at (170bp,152bp) [place] {$J_7$};
  \node (6) at (150bp,220bp) [place] {$J_6$};
  \node (9) at (50bp,150bp) [place] {$J_9$};
  \node (8) at (271bp,139bp) [place] {$J_8$};
  \draw [thin] (5) -- (8);
  \draw [thin] (6) -- (7);
  \draw [dashed,thin] (8) -- (10);
  \draw [thin] (3) -- (4);
  \draw [thin] (9) -- (11);
  \draw [dashed,thin] (7) -- (11);
  \draw [thin] (2) -- (3);
  \draw [thin] (8) -- (11);
  \draw [thin] (9) -- (10);
  \draw [thin] (7) -- (10);
  \draw [thin] (6) -- (8);
  \draw [thin] (4) -- (6);
  \draw [thin] (4) -- (5);
  \draw [thin] (11) -- (12);
  \draw [thin] (6) -- (9);
  \draw [thin] (4) -- (8);
  \draw [thin] (7) -- (8);
  \draw [thin] (4) -- (7);
  \draw [thin] (3) -- (5);
  \draw [thin] (1) -- (2);
  \draw [thin] (10) -- (11);
  \draw [thin] (1) -- (3);
  \draw [loosely dotted,thick] (20bp,30bp) circle (40pt);
  \draw [loosely dotted,thick,rounded corners=50pt] (-10bp,200bp) -- (290bp,75bp) -- (80bp,42bp) -- cycle;
  \draw [loosely dotted,thick,rounded corners=20pt] (140bp,133bp) -- (115bp,240bp) -- (425bp,360bp) -- (600bp,300bp) -- (460bp,170bp) -- (300bp,115bp) -- cycle;
  \draw [very thick,black!25,rounded corners=20pt] (65bp,65bp) -- (150bp,180bp) -- (315bp,155bp) -- (240bp,60bp) -- cycle;
 \end{tikzpicture}
\end{center}
\caption{\label{fig:BorelIncidence} The Borel incidence graph of $\hilb_{4t+1}^4$. The dashed lines correspond to composed rational deformations. With the dotted lines, we enclose the Borel-fixed ideals lying on a common component because of Reeves criterion \cite{RA}, whereas with the thick gray line we enclose the Borel ideals lying on a common component because of our criterion.} 
\end{figure}
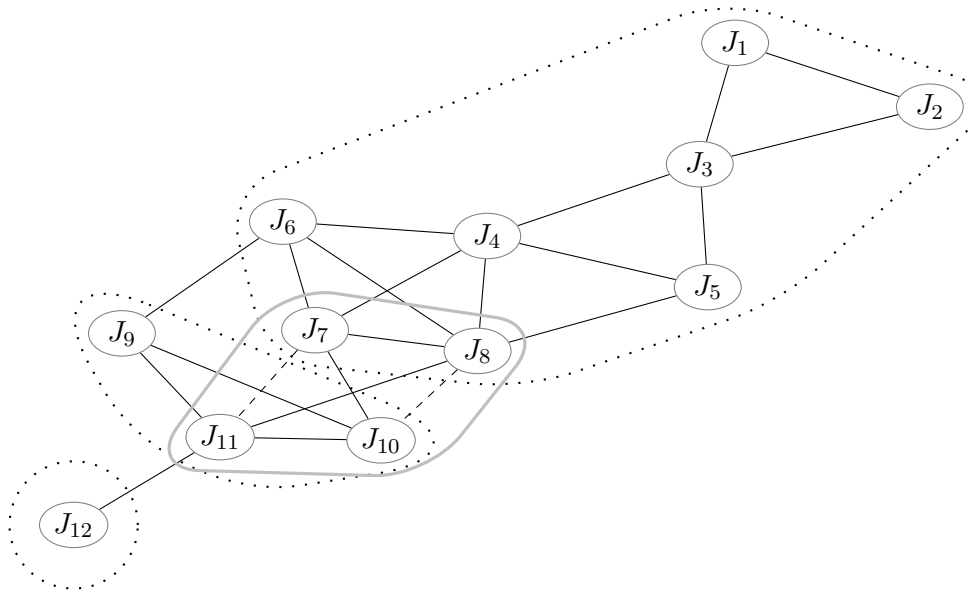
\end{example}


\providecommand{\bysame}{\leavevmode\hbox to3em{\hrulefill}\thinspace}
\providecommand{\MR}{\relax\ifhmode\unskip\space\fi MR }
\providecommand{\MRhref}[2]{%
  \href{http://www.ams.org/mathscinet-getitem?mr=#1}{#2}
}
\providecommand{\href}[2]{#2}

\end{document}